\documentclass{amsart}

\usepackage{latexsym}
\usepackage{amsfonts,amssymb}
\usepackage{pstricks,pst-plot}
\usepackage{color}
\usepackage{pstcol,pst-node}

\newtheorem{theorem}{Theorem}[section]
\newtheorem{lemma}[theorem]{Lemma}
\newtheorem{corollary}[theorem]{Corollary}

\newtheorem{remark}[theorem]{Remark}
\theoremstyle{definition}
\newtheorem{definition}[theorem]{Definition}
\newtheorem{example}[theorem]{Example}

\theoremstyle{remark}

\begin{document}

\newcommand{\abs}[1]{\lvert#1\rvert}
\newcommand{\ID}{\mathbb{I}}
\newcommand{\complessi}{\mathbb C}
\newcommand {\interi}{\mathbb Z}
\newcommand{\reali}{I\!\!R}
\newcommand{\razionali}{\mathbb Q}
\newcommand{\pd}{\partial_x}
\newcommand\Pa{Painlev\'e}
\newcommand\Pt{{\rm P}_{\rm{\scriptstyle I}}}
\newcommand\Lr{{\mathcal L}}
\newcommand{\nn}{\nonumber}
\newcommand\Pol{{\mathcal P}}
\newcommand\Se{{\mathcal S}}
\newcommand\A{{\mathcal A}}
\newcommand\M{{\mathsf M}}
\newcommand {\KZ}{Knizhnik--Zamolodchikov }
\newcommand {\KZt}{Knizhnik--Zamolodchikov--type }
\newcommand {\KKZ}{_{\scriptscriptstyle{\operatorname{KZ}}}}
\newcommand {\nablak}{\nabla\KKZ}
\def\ddt{\frac{\rm d}{{\rm d} t}}
\def\tr{{\rm Tr}}
\def\res{{\rm res}}
\def\ddx{\frac{\rm d}{{\rm d} x}}
\def\ddtt{\frac{\rm d}{{\rm d}\tau}}
\def\ddz{\frac{\rm d}{{\rm d} z}}
\def\RE{{\rm Re}}
\def\IM{{\rm Im}}
 \def\arg{{\rm arg}}
\def\iff{\Longleftrightarrow}
\def\U{{\mathcal U}}
\def\pderiv{\partial_x}
\def\orexp#1{#1^{\raise3pt\hbox{$\scriptstyle\!\!\!\!\!\!\!\to$}}\,{}}
\def\one#1{#1^{\raise5pt\hbox{$\scriptstyle\!\!\!\!1$}}\,{}}
\def\two#1{#1^{\raise5pt\hbox{$\scriptstyle\!\!\!\!2$}}\,{}}
\def\onetwo#1{#1^{\raise4pt\hbox{$\scriptstyle\!\!\!\!\!{12}$}}\,{}}
\def\twoone#1{#1^{\raise4pt\hbox{$\scriptstyle\!\!\!\!\!{21}$}}\,{}}

\def\GL{\mathcal{GL}}

\title[Quantum Teichm\"uller theory]
{Lecture Notes on Quantum Teichm\"uller  theory.}

\author{Leonid Chekhov.}

\maketitle

\begin{center}{\it Notes collected by M. Mazzocco}\end{center}

\begin{abstract}  These notes are based on a lecture course by L. Chekhov held at the
University of Manchester in May 2006 and February-March 2007. They
are divulgative in character, and instead of  containing rigorous
mathematical proofs, they illustrate statements giving an intuitive
insight. We intentionally remove most bibliographic references from the body
of the text devoting a special section to the history of the subject
at the end.
\end{abstract}

\tableofcontents

\section{Combinatorial description of Teichm\"uller spaces}

In this lecture we concentrate on the Teichm\"uller space $\mathcal
T_{g}^s$ of  Riemann surfaces $F_{g,s}$ of genus $g$ with $s$ holes.

\begin{figure}[h]
\begin{center}\psset{unit=2.5cm}
\begin{pspicture}(-1.5,-0.4)(3,0.5)
 \parametricplot{-1.16}{1.16}{t 2 exp 1.5 t mul t 3 exp -2.2 mul add t 5 exp add}
 \psellipse(1.42,0)(0.2,0.44) 
 \psarc(0.5,0.3){.4}{230}{310}
 \psarc(0.5,-0.4){.4}{60}{120}
\end{pspicture}
  \caption{A torus with one hole}\label{torus-hole}
\end{center}
\end{figure}
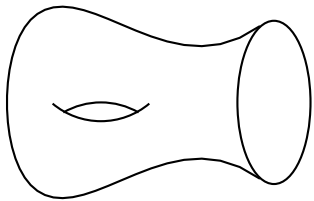

\begin{example}
The simplest case is the torus with one hole, $F_{1,1}$ (see figure \ref{torus-hole}).
 \end{example}

The holes can be considered in different parameterizations. For this
lecture course we concentrate on the Poincar\'e uniformization. There is
another uniformization which is due to Strebel, but in the
Poincar\'e one there is a good description of the Poisson structure
and of the quantization.

What is the Poincar\'e unformization? On each $F_{g,s}$ we may
define several metrics. We consider two metrics to be equivalent if
they are mapped one to another by a diffeomorphism. Within each
equivalence class, we chose the representative to be a metric with
local constant curvature $-1$.  Observe that for $g>1$ there always
exists such representative element (i.e. we may pick $s=0$). For
$g=1$ we need at least one hole, otherwise it is impossible as the
universal cover of the torus is flat. For $g=0$, we need at least
$s=3$.

This is called {\it Poincar\'e uniformization:}\/ our Riemann
surface $F_{r,s}$ is mapped to a surface with local constant
curvature $-1$, namely
$$
F_{g,s} = {\mathbb H}\slash\Delta_{g,s},
$$
where $\mathbb H$ denotes the upper half plane and $\Delta_{g,s}$ is
a {\it Fuchsian group,}\/ i.e. a finitely generated discrete
subgroup of the isometry group $\mathbb PSL(2,\mathbb R)$ of
$\mathbb H$:
$$
\Delta_{g,s} =\langle\gamma_1,\dots,\gamma_{2g+s}\rangle,
$$
where $\gamma_1,\dots,\gamma_{2g+s}$ are {\it hyperbolic
elements,}\/ i.e. they have two distinct fixed points on the {\it
absolute,}\/ i.e. $\mathbb R\cup\{\infty\}$.

Recall that on $\mathbb H$ we have the metric ${\rm d}s^2=\frac{{\rm
d} x^2+{\rm d} y^2}{y^2}$ and the geodesics are either semi--circles
with centre on the $x$--axis or half--lines parallel to the
$y$--axis (see figure \ref{pic:geodesics}). To familiarize the
reader with the hyperbolic metric, we have discussed some examples
in Appendix~A.
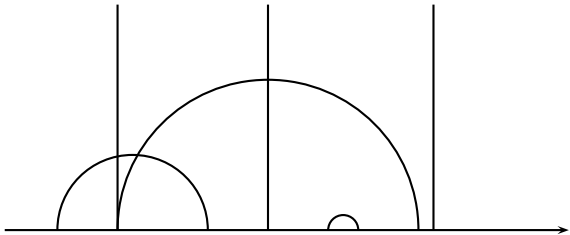
\begin{figure}[h]
\begin{center}
\begin{pspicture}(-2,0)(2,3)
 \psline{->}(-3.5,0)(4,0)
 \psline(-2,0)(-2,3)
 \psline(0,0)(0,3)
 \psline(2.2,0)(2.2,3)
 \psarc(-1.8,0){1}{0}{180}
 \psarc(0,0){2}{0}{180}
 \psarc(1,0){.2}{0}{180}
\end{pspicture}
\caption{Some geodesics in the upper half plane}\label{pic:geodesics}
\end{center}
\end{figure}

Observe that given a geodesic $\gamma$ and a point $P$ in $\mathbb H$ not
belonging to it, there are infinitely many geodesics through $P$
which do not cross $\gamma$ (therefore Euclid's postulate does not hold
true in hyperbolic geometry). The $x$--axis is infinitely distant
from the points in $\mathbb H$, this is why we say that it belongs
to the {\it absolute.}

Note that actually, no ``points'' lie on the absolute; instead we consider classes of
geodesics terminating at such a point in a sense that all of them become asymptotically
close at large proper distances: the collection of such ``points'' is the open real line together
with the infinity point (defined by upper ends of vertical geodesic lines: all these lines in
Fig.~\ref{pic:geodesics} are asymptotically close). The equivalent and often used picture is the
Poincar\'e disc: the absolute is the boundary of the disk and the geodesics are arcs.

Let $A\in\mathbb P SL(2,\mathbb R)$ with $A=\begin{pmatrix} a&b\cr c&d\end{pmatrix}$.
This acts on $\mathbb H$ by the M\"obius transformation
$$
\gamma_A : z \mapsto \frac{az+b}{cz+d}.
$$
This action is {\it transitive,}\/ and  the M\"obius transformation has two fixed points (possibly at infinity),
given by
$$
z_\pm = \frac{1}{2c}\left(a-d \pm \sqrt{T^2-4}\right),
$$
where $T=\mathrm{trace}(A)$. If one fixed point is at $\infty$
(i.e.\ if $c=0$) then the other is at $b/(d-a)\in\mathbb R$. Notice
that $T$ is not well-defined on $\mathbb P SL(2,\mathbb R)$ but
$T^2$ is. Analogously, the ratio $\frac{a-d}{c}$ is well defined in
$\mathbb P SL(2,\mathbb R)$, so the fixed points are indeed uniquely
determined by the element in $\mathbb P SL(2,\mathbb R)$. Two
non--identity elements in  $\mathbb P SL(2,\mathbb R)$ commute if
and only if they have the same fixed points.

The element $A\in \mathbb P SL(2,\mathbb R)$ is said to be \emph{hyperbolic} if
$T^2>4$, \emph{elliptic} if $T^2<4$ and \emph{parabolic} if $T^2=4$.
Parabolic elements have a unique fixed point on the absolute, and are conjugate to
$\begin{pmatrix}1&1\cr 0&1\end{pmatrix}$.

We will almost exclusively be interested in hyperbolic elements. The
fixed points of such elements are real (lie on the absolute). Since
the M\"obius transformation $\gamma_A$ is a hyperbolic isometry, it
follows that it preserves the unique geodesic between its two fixed
points. We call such geodesic {\it invariant axis.}

The eigenvalues of the matrix $A$ are given by
$$\lambda_\pm = \frac12 \left( T \pm  \sqrt{T^2-4}\right).$$
The linearization of the M\"obius transformation $\gamma_A$, at the
fixed point $z_\pm$ is the complex linear map with eigenvalue
$\lambda_\pm^2$.

Given any element $\gamma\in\mathbb P SL(2,\mathbb R)$, $\gamma=\left(\begin{array}{cc} a&b\\ c&d
\end{array}\right)$, $ad-bc=1$, we can uniquely determine it by its eigenvalues and its eigenvectors.
We are now going to characterize $\gamma$ by two other objects: a closed geodesic and its length.

In fact, since the determinant of $\gamma$ is one, the eigenvalues
can be expressed as $\exp\left(\pm\frac{l_\gamma}{2}\right)$. In the
basis of the eigenvectors, $\gamma(z)=\exp(l_\gamma) z$, so $\gamma$
is a dilation and because it must be an isometry, it maps the whole
geodesic though $z$ to a second geodesic (see figure \ref{pic:hyp-action}).

\begin{figure}[h]
\begin{center}
\begin{pspicture}(-3.5,0)(9,3)
 \psline{->}(-3.5,0)(-0.2,0)
 \psarc(-1.8,0){1}{0}{180}
 \psline{->}(1,0)(9,0)
 \pscustom{
   \psarcn(4.8,0){3}{180}{0}
 }
 \psline{|->}(0.2,0.5)(0.8,0.5)
\end{pspicture}
\caption{}\label{pic:hyp-action}
\end{center}
\end{figure}

Our diagonalized element $\gamma$ has two distinct fixed points which are $0$ and
$\infty$. Generally if $\gamma$ is not diagonal, the position of the
fixed points is uniquely determined by the eigenvectors. They always
lie on the absolute.

\begin{example}\label{ex:hypnd}
Consider $\gamma=\left(\begin{array}{cc}2&2\\ \frac{1}{2}&1\end{array}\right)$,
it has eigenvalues $\lambda_\pm = \frac{1}{2}\left(3\pm\sqrt{5}\right)$ and eigenvectors
$v_\pm = (1\mp\sqrt{5},1)$. So the fixed points lie at $z=1\mp\sqrt{5}$.
The other geodesics are mapped as in figure \ref{pic:hyp-action1}.
\begin{figure}[h] 
\begin{center}
\psset{unit=8mm}
\begin{pspicture}(6,-1)(9,4)
 \psline{->}(0.5,0)(6.4,0)
  \psarcn[linecolor=blue,linewidth=2pt](3.2,0){2}{180}{0}
 \psline[linecolor=green,linewidth=2pt](5.9,0)(5.9,4)
\rput(7,0.5){$\leadsto$}
  \psline{->}(8,0)(15,0)
  \psarc[linecolor=green,linewidth=2pt](10,0){1.5}{0}{180}
 \psarc[linecolor=blue,linewidth=2pt](13.5,0){1}{0}{180}
\end{pspicture}\caption{}\label{pic:hyp-action1}
\end{center}
\end{figure}
\end{example}

If we identify points on  $\mathbb H$ by the action of $\gamma$, it
means that we have to identify the initial geodesic and its image
under $\gamma$, so that we obtain an infinite hyperboloid.  Let us
consider the only finite (i.e. not connecting points which lie on
the absolute) geodesic contained between the initial geodesic and
its $\gamma$--image. This is mapped to the only closed geodesic in
our hyperboloid (figure \ref{pic:closed-geo}).

\begin{figure}[h] 
\begin{center}\psset{unit=8mm}
\begin{pspicture}(1,-2)(9,4)
 \psline{->}(1,0)(9,0)
 \psarc(4.8,0){1}{0}{180}
 \psline(3.8,0)(1.8,0)
 \psarcn(4.8,0){3}{180}{0}
 \psline(7.8,0)(5.8,0)
 \psline[linecolor=red,linewidth=2pt](4.8,1)(4.8,3)
   \rput(10,0.5){$\leadsto$}
\end{pspicture}
\begin{pspicture}(-4,-3)(2,4)
 \parametricplot{-1}{1}{1 t 2 exp add t 3 mul}   
 \parametricplot{-1}{1}{-1 1 t 2 exp add mul t 3 mul} 
 \parametricplot[linecolor=red,linestyle=dashed,linewidth=2pt]{0}{180}{t cos 0.3 t sin mul}
 \parametricplot[linecolor=red,linewidth=2pt]{180}{360}{t cos 0.3 t sin mul}
 \psellipse(0,3)(2,0.6)
 \parametricplot[linestyle=dashed]{0}{180}{2 t cos mul 0.6 t sin mul -3 add}
 \parametricplot{180}{360}{2 t cos mul 0.6 t sin mul -3 add}
\end{pspicture}
\caption{}\label{pic:closed-geo}
\end{center}
\end{figure}

The length of this closed geodesic is
$$
\int_1^{l_\gamma} \frac{{\rm d}y}{y}=l_\gamma,
$$
where we have assumed for sake of simplicity that the initial
geodesic had Euclidean radius one. So we see that for each
hyperbolic element $\gamma$ there exists a unique closed geodesic
and its hyperbolic length is related to the trace of $\gamma$:
$$
{\rm Tr}(\gamma)=2 \cosh l_\gamma.
$$
As the trace is invariant under the action of $\mathbb P
SL(2,\mathbb R)$, the hyperbolic--length is a {\it
hyperbolic--invariant.}\/ Functions on $\mathbb H$ which are
invariant under the action of a Fuchsian group are called {\it
automorphic functions.}\/ Observe that the red closed geodesic lies
on a geodesic of $\mathbb H$ which relates two points: $0$ and
$\infty$. These are exactly the stable points of $\gamma$.

\begin{lemma}
Closed geodesics are in one-to-one correspondence with invariant
axis of hyperbolic translations.
\end{lemma}

Observe that actually we are interested in the Teichm\"uller space.
Therefore the whole construction is invariant under the action of
the automorphism group of $\mathbb H$, that is,  $\mathbb P SL(2,\mathbb R)$.
The only relevant information about $\gamma$ is that it has two
distinct eigenvalues, and its trace. So we can state the following

\begin{lemma}\label{1to1.1}
The following objects are in one--to--one correspondence
\begin{enumerate}
\item Conjugacy classes of hyperbolic elements $\gamma\in\Delta_{g,s}$.
\item Closed geodesics on a Riemann surface; their lengths are related to
the traces of $\gamma$.
\end{enumerate}
\end{lemma}

In virtue of the above lemma we shall often identify conjugacy
classes of hyperbolic elements with their representative $\gamma$
and we shall denote the associated closed geodesic in the same way,
i.e. $\gamma$.

We now show that they are in one to one correspondence with the
homotopy classes in $\pi_1(F_{g,s})$. In fact in each homotopy class
we have only one shortest representative. Its length determines
$l_\gamma$. Its topology determines the invariant points of
$\gamma$.

\begin{figure}[h]
\begin{center}
\begin{pspicture}(4,0)(12,4)
 \psline{->}(4,0)(12,0)
 \psarc(8,0){1}{0}{180}
 \psarc(8,0){3}{0}{180}
 \psarc[linecolor=green,linewidth=2pt](10,0){0.5}{0}{180}
 \psarc[linecolor=yellow,linewidth=2pt](6,0){0.5}{0}{180}
 \psarc[linecolor=blue,linewidth=2pt](8,0){2}{15}{165}
 \psline[linecolor=red,linewidth=2pt](8,1)(8,3)
\end{pspicture}
\end{center}
 \caption{The red invariant axis corresponds to $\gamma$, while the blue one to $\tilde \gamma$.
}\label{pic:cl-geo}
\end{figure}
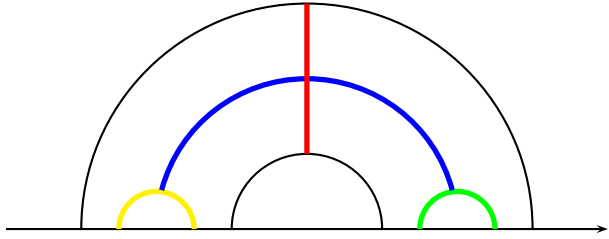

To show this, we illustrate how to reconstruct the Riemann surface
$F_{g,s}$ from the generators $\gamma_1,\dots,\gamma_n$ of the
Fuchsian group $\Delta_{g,s}$. We do this in the simplest possible case,
i.e. when there are only two hyperbolic generators, $\gamma$ and
$\tilde\gamma$. We assume that $\gamma$ has fixed points $0$ and
$\infty$ as above and that $\tilde\gamma$ has fixed points on the
the $x$--axis. In both cases the geodesic connecting the fixed
points is called {\it invariant axis}. The action of $\tilde\gamma$
on the geodesics is supposed to be like in example \ref{ex:hypnd}.
To understand better let us draw the picture on the Poincar\'e disk.
In this setting each arc is mapped in the opposite one and
its {\em inside} (white domain) is mapped to the {\em outside} of the opposite one. In the figure
\ref{pic:poincare0} the green arc is mapped to the yellow one by
$\tilde\gamma$. We draw also the same geodesics in the hyperbolic
plane. The shaded area is the {\it fundamental domain.} In figure \ref{pic:cl-geo} we show the invariant axis.

\begin{figure}[h]
\begin{center}
\begin{pspicture}(-1,-1)(2,4)

\pscircle[fillstyle=crosshatch] (0,1.5){2}

\psarc[linestyle=none,fillstyle=solid,fillcolor=white](2.5,1.5){1.5}{100}{250}
\psarc[linestyle=none,fillstyle=solid,fillcolor=white](-2.5,1.5){1.5}{-70}{70}
\psarc[linestyle=none,fillstyle=solid,fillcolor=white](0,4.3){1.96}{200}{330}
\psarc[linestyle=none,fillstyle=solid,fillcolor=white](0,-1){1.5}{20}{160}
\pscircle(0,1.5){2}
 \psarc[linecolor=green,linewidth=2pt](2.5,1.5){1.5}{127}{233}
 \psarc[linecolor=yellow,linewidth=2pt](-2.5,1.5){1.5}{-53}{53}
 \psarc(0,4.3){1.96}{224.4}{315.6}
 \psarc(0,-1){1.5}{37}{143}
 \end{pspicture}
  \begin{pspicture}(3.7,-1)(11.5,3.5)
 \psline{->}(4.5,0)(11.5,0)
 \psarc[fillstyle=crosshatch](8,0){3}{0}{180}
  \psarc[fillstyle=solid,fillcolor=white](8,0){1}{0}{180}
 \psarc[linecolor=green,linewidth=2pt,fillstyle=solid,fillcolor=white](10,0){0.5}{0}{180}
 \psarc[linecolor=yellow,linewidth=2pt,fillstyle=solid,fillcolor=white](6,0){0.5}{0}{180}
\end{pspicture}
\end{center}
 \caption{}\label{pic:poincare0}
\end{figure}

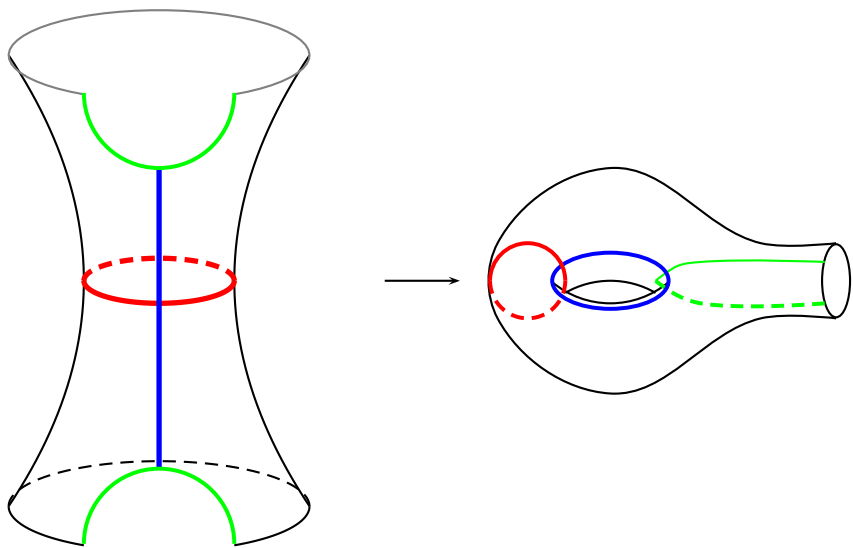
\begin{figure}[h] 
\begin{center}
\begin{pspicture}(-2,-4)(10,4)
 \parametricplot{-1}{1}{1 t 2 exp add t 3 mul}   
 \parametricplot{-1}{1}{-1 1 t 2 exp add mul t 3 mul} 
 \parametricplot[linecolor=red,linestyle=dashed,linewidth=2pt]{0}{180}{t cos 0.3 t sin mul}
 \parametricplot[linecolor=red,linewidth=2pt]{180}{360}{t cos 0.3 t sin mul}
 \psline[linecolor=blue,linewidth=2pt](0,-2.5)(0,1.5)
 \parametricplot[linecolor=gray]{-60}{240}{t cos 2 mul 0.6 t sin mul 3 add}
 \parametricplot[linestyle=dashed]{0}{180}{2 t cos mul 0.6 t sin mul -3 add}
 \parametricplot{180}{240}{2 t cos mul 0.6 t sin mul -3 add}
 \parametricplot{-60}{0}{2 t cos mul 0.6 t sin mul -3 add}
 \psarc[linecolor=green,linewidth=1.5pt](0,2.5){1}{-180}{0}
 \psarc[linecolor=green,linewidth=1.5pt](0,-3.5){1}{0}{180}
 \psline{->}(3,0)(4,0)
 \pscurve(9,-0.5)(8,-0.5)(6,-1.5)(4.5,-0.5)(4.5,0.5)(6,1.5)(8,0.5)(9,0.5)
 \psellipse(9,0)(0.2,0.5)
 \psarc(6,0.9){1.2}{230}{310}
 \psarc(6,-1.2){1.2}{60}{120}
 \pscurve[linecolor=green](8.85,0.25)(7.2,0.25)(6.9,0.2)(6.6,0)
 \pscurve[linestyle=dashed,linecolor=green,linewidth=1.5pt](8.85,-0.3)(7.2,-0.3)(6.9,-0.2)(6.6,0)
 \psellipse[linecolor=blue,linewidth=1.5pt](6,0)(0.8,0.4)
 \psarc[linecolor=red,linewidth=1.5pt](4.9,0){0.5}{0}{180}
 \psarc[linecolor=red,linestyle=dashed,linewidth=1.5pt](4.9,0){0.5}{180}{360}
\end{pspicture}
\end{center}
 \caption{Glue along the green curve}\label{pic:hyp-to}
\end{figure}

We now glue as in figure \ref{pic:closed-geo}, obtaining a
hyperboloid. The yellow and green geodesic are drawn on the
hyperboloid in figure \ref{pic:hyp-to}. Since we are going to identify them, we draw them in
the same colour (green). By gluing along the green geodesics we get
our torus with one hole (see figure \ref{pic:hyp-to}).
Observe that in particular we can identify the red geodesic with
an $A$-cycle and the blue one with a $B$-cycle. To remember this we shall denote
$$
\gamma_A = \gamma,\qquad
\gamma_B=\tilde\gamma.
$$

Once we have obtained our torus with one hole, it is clear that only
the homotopy class of the closed geodesics corresponding  to the
invariant axes (i.e. the red and blue ones) are relevant.

\begin{lemma}\label{1to1}
The following objects are in one--to--one correspondence:
\begin{enumerate}
\item Conjugacy classes of hyperbolic elements $\gamma\in\Delta_{g,s}$.
\item Closed geodesics of given length.
\item Conjugacy classes in $\pi_1(F_{g,s})$.
\end{enumerate}
\end{lemma}

Finally we want to prove that all this is in one-to-one with the so
called {\it fat--graph.}\/ This is going to give us a combinatorial
description of the Teichm\"uller space.

\begin{figure}[h] 
\begin{center}
\begin{pspicture}(-1.5,-2)(2.5,2) 
 \psline(0,-2)(0,2)
 \psline(-1.5,0)(1.5,0)
\psline{->}(2,0)(2.5,0)
\end{pspicture}
\begin{pspicture}(-1.8,-2)(2.5,2) 
 \psline(0,-2)(0,2)
 \psline(-1.5,0.4)(0,0.4)
 \psline(0,-0.4)(1.5,-0.4)
 \psline{->}(2,0)(2.5,0)
\end{pspicture}
\begin{pspicture}(-2,-2)(2,2) 
\newgray{mygray}{0.96}
{ \psset{arrowsize=6pt,linewidth=1.5pt}
  \psline(-0.2,2)(-0.2,0.6)(-1.7,0.6)
  \psline(-1.7,0.2)(-0.2,0.2)(-0.2,-2)
  \psline(0.2,-2)(0.2,-0.6)(1.7,-0.6)
  \psline(1.7,-0.2)(0.2,-0.2)(0.2,2)
 \psline{->}(-0.2,0.6)(-1,0.6)\psline{->}(-0.2,0.2)(-0.2,-1)\psline{->}(0.2,-0.6)(1,-0.6)\psline{->}(0.2,-0.2)(0.2,1)
} \psarcn{->}(-1.3,-1){0.5}{120}{-30}
 \psarcn{->}(-1.3,1.5){0.5}{20}{-120}
 \psarcn{->}(1.3,1){0.5}{-30}{-210}
 \psarcn{->}(1.3,-1.5){0.5}{210}{30}
\end{pspicture}
\end{center}
 \caption{Resolving of a crossing into a fat graph}\label{pic:ribbon}
\end{figure}
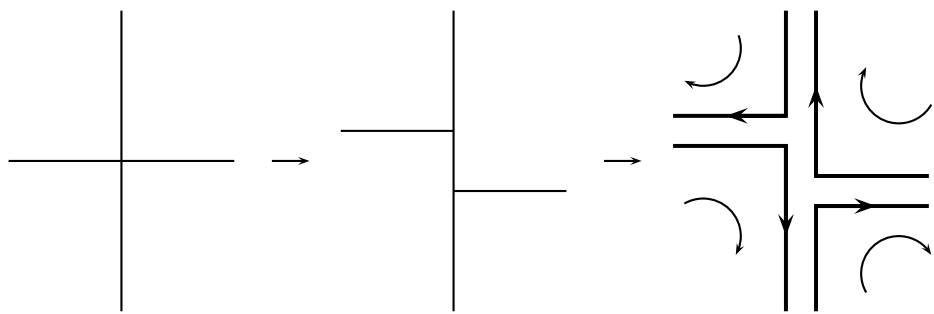

Usually graphs are determined by their incidence matrix specifying
how many edges meet at each vertex. Suppose the faces carry an
orientation. In order to be able to put an orientation on the edges
which is compatible with the orientation of the faces, we make the
graph fat, or ribbon. Moreover we shift each crossing so
that at each vertex only 3 edges come together (see figure
\ref{pic:ribbon}).

In this way, a curve coming from an edge at a vertex may turn either
left or right (which makes all computations combinatorially
simpler).

\begin{figure}[h] 
\begin{center}
\begin{pspicture}(4,-1.5)(9,2)
 \psecurve(9,-0.5)(8,-0.5)(6,-1.5)(4.5,-0.5)(4.5,0.5)(6,1.5)(8,0.5)(9,0.5)
 \psellipse(8,0)(0.2,0.5)
 \psarc(6,0.9){1.2}{230}{310}
 \psarc(6,-1.2){1.2}{60}{120}
 \parametricplot[linecolor=red,linewidth=1.5pt]{0}{180}{t sin -0.2 mul 6 add t cos 0.75 mul 0.75 add}
 \parametricplot[linecolor=red,linewidth=1.5pt,linestyle=dashed]{0}{180}{t sin 0.2 mul 6 add t cos 0.75 mul 0.75 add}
 \psellipse[linecolor=blue,linewidth=1.5pt](6,0)(0.8,0.4)
 \psline{->}(8.5,0)(9,0)
\end{pspicture}
\begin{pspicture}(5,-1.5)(7,2)
 \psellipse[linecolor=red,linewidth=1.5pt](6,0.75)(0.2,0.75)
 \pscircle[fillstyle=solid,fillcolor=white,linecolor=white](6.2,0.35){0.15}
 \psellipse[linecolor=blue,linewidth=1.5pt](6,0)(0.8,0.4)
 \psline{->}(7,0)(7.5,0)
\end{pspicture}
\begin{pspicture}(4,-1.5)(7,2)
 \pscurve[linecolor=red,linewidth=1.5pt,showpoints=false](5.2,.4)(5.1,1.2)(6,1.3)(6.3,0.3)(7,0)
 \pscircle[fillstyle=solid,fillcolor=white,linecolor=white](6.2,0.75){0.15}
 \psellipse[linecolor=blue,linewidth=1.5pt](6,0)(1,0.8)
 \psline{->}(7.5,0)(8,0)
\end{pspicture}
\begin{pspicture}(3.5,-1.5)(7,2)
 \pscurve[linecolor=red,linewidth=1.5pt,showpoints=false](5.2,.6)(5.05,1.45)(6,1.6)(6.2,0.2)(6.9,0.3)
 \pscurve[showpoints=false](5,.6)(4.9,1.6)(6.15,1.75)(6.35,0.35)(6.7,0.4)
 \pscurve[showpoints=false](5.3,.95)(5.2,1.4)(5.8,1.5)(6,0.1)(6.8,0.1)
 \pscircle[fillstyle=solid,fillcolor=white,linecolor=white](6.1,1){0.35}
 \psarc(6,0){1.2}{150}{125}
 \psarc(6,0){0.8}{30}{5} %
 \pscircle[linecolor=blue,linewidth=1.5pt](6,0){1}
\end{pspicture}
\end{center}
 \caption{How to associate a fat-graph to a Riemann surface.}\label{pic:torus-to-pr}
\end{figure}
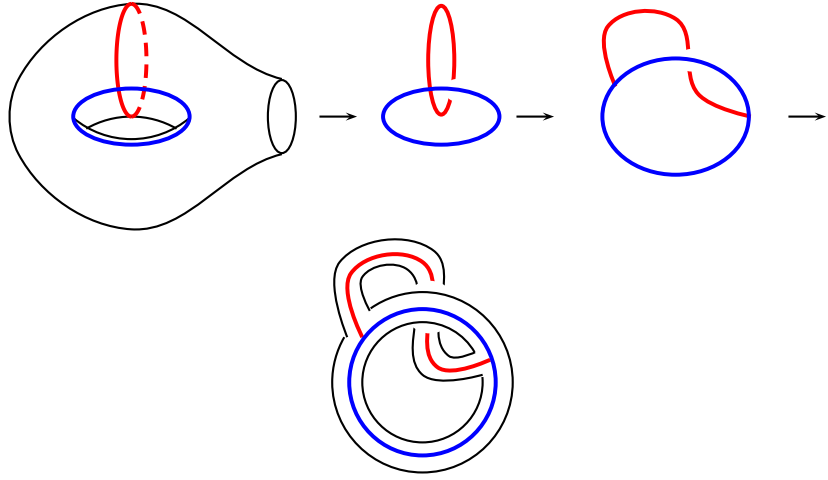

Now in our torus with one hole let us pick the graph defined by our $A$ and $B$
cycles and make it into a fat graph. Observe that the this fat graph
has only one face and such face contains the hole (see figure \ref{pic:torus-to-pr}).

So the rules for
a fat graph are
\begin{enumerate}
\item The vertices are trivalent
\item Each face must contain only one hole.
\end{enumerate}

For each such fat graph $\Gamma$ there exists a unique Riemann
surface for which $\Gamma$ gives the minimal simplicial
decomposition. Then we can add a fourth element in lemma \ref{1to1}:

\begin{theorem}\label{1to1-final}
The following objects are in one--to--one correspondence:
\begin{enumerate}
\item Conjugacy classes of hyperbolic elements $\gamma\in\Delta_{g,s}$.
\item Closed geodesics of given length.
\item Conjugacy classes in $\pi_1(F_{g,s})$.
\item Closed paths in a fat graph $\Gamma$.
\end{enumerate}
\end{theorem}

\section{Coordinates in the Teichm\"uller space}

We now want to find good coordinates in the Teichm\"uller space
${\mathcal T}_g^s$. Recall that each point in ${\mathcal T}_g^s$ is
specified by a Fuchsian group $\Delta_{g,s}$ generated by hyperbolic elements each
of which is uniquely determined by its invariant axis in the
fundamental domain. Given a hyperbolic element $\gamma$, we are
going to express it as the product of matrices of type
$$
L:=\left(\begin{array}{cc}0&1\\-1&-1\\
\end{array}\right), \quad
R:=\left(\begin{array}{cc}1&1\\-1&0\\
\end{array}\right), \quad
X_{z}:=\left(\begin{array}{cc}0&-\exp\left({\frac{z}{2}}\right)\\
\exp\left(-{\frac{z}{2}}\right)&0\\
\end{array}\right),
$$
for some $z\in{\mathbb R}$. Projectively $R=L^2$ and $L=R^2$. We are
looking for a decomposition of the form
$$
\gamma = R^{k_n} X_{z_n} \dots R^{k_1} X_{z_1},
$$
where $k_i$ can be either $1$ or $2$. Observe that $\gamma$ obtained
in this way will necessarily be real and have the following sign
structure
$$
{\rm sign}(\gamma)=\left(\begin{array}{cc}
+&-\\
-&+\\
\end{array}\right)
$$
because each of the matrices $RX_{z_i}$ and $LX_{z_i}$ has this
structure. This is not a restrictive condition, because we are
taking $\gamma\in{\mathbb P}SL(2,\mathbb R)$ with freedom of
conjugation by ${\mathbb P}SL(2,\mathbb R)$ as due to Theorem
\ref{1to1-final}, only the conjugacy class of $\gamma$ is relevant.

\begin{example}\label{ex:gammab}
Let us factorize an element of the Fuchsian group $\gamma$ in the above way. Suppose $\gamma$ maps
$$
0\to -1,\qquad \exp(z_0) \to \infty,\qquad
\infty\to -1-\exp(-z_1),
$$
for some real numbers $z_0,z_1$.  It is a straightforward computation to show that
$$
\gamma =R X_{z_1} L X_{z_0},
$$
is the correct factorization.
\end{example}

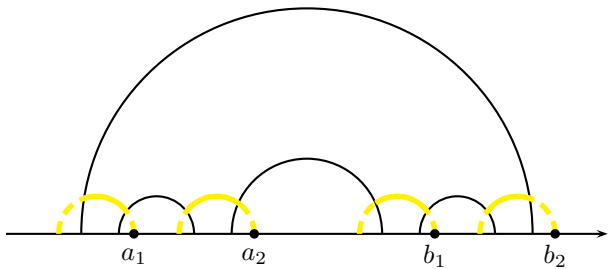
\begin{figure}[h]
\begin{center}
\begin{pspicture}(4,-0.5)(12,4)
 \psline{->}(4,0)(12,0)
 \psarc(8,0){1}{0}{180}
 \psarc(8,0){3}{0}{180}
 \psarc(10,0){0.5}{0}{180}
 \psarc(6,0){0.5}{0}{180}
 \psarc[linecolor=yellow,linewidth=2pt](5.2,0){0.5}{33}{110}
  \psarc[linecolor=yellow,linewidth=2pt,linestyle=dashed](5.2,0){0.5}{0}{180}
\psdot(5.7,0)  \rput(5.7,-0.3){$a_1$}
 \psarc[linecolor=yellow,linewidth=2pt](6.8,0){0.5}{50}{145}
\psarc[linecolor=yellow,linewidth=2pt,linestyle=dashed](6.8,0){0.5}{0}{180}
\psdot(7.3,0)  \rput(7.3,-0.3){$a_2$}
\psarc[linecolor=yellow,linewidth=2pt](9.2,0){0.5}{40}{125}
\psarc[linecolor=yellow,linewidth=2pt,linestyle=dashed](9.2,0){0.5}{0}{180}
\psdot(9.7,0)  \rput(9.7,-0.3){$b_1$}
 \psarc[linecolor=yellow,linewidth=2pt](10.8,0){0.5}{65}{145}
 \psarc[linecolor=yellow,linewidth=2pt,linestyle=dashed](10.8,0){0.5}{0}{180}
 \psdot(11.3,0)  \rput(11.3,-0.3){$b_2$}
\end{pspicture}
\caption{Image of the bottleneck curve in the fundamental domain.}\label{pic:bottleneck}
\end{center}
\end{figure}

We are now going to show a way to guess such factorization by
constructing an {\it ideal triangulation}\/ of the fundamental
domain of the Riemann surface in which all the infinite hyperboloids
attached to the holes have been cut off. More precisely, consider
the closed curves homeomorphic to the holes, say the {\it bottleneck
curves} on each hyperboloid. These are closed geodesics, therefore
they correspond to an element of the Fuchsian group. In the case of
one hole in the torus the image of the bottleneck curve in the
fundamental domain is given in figure \ref{pic:bottleneck}.

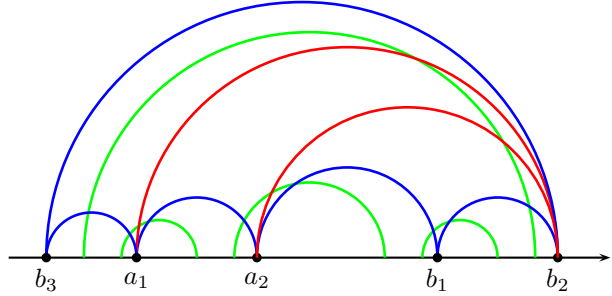
\begin{figure}[h]
\begin{center}
\begin{pspicture}(1,-1)(16,4)
 \psline{->}(4,0)(12,0)
  \psarc[linecolor=green,linewidth=1pt](8,0){1}{0}{180}
 \psarc[linecolor=green,linewidth=1pt](8,0){3}{0}{180}
 \psarc[linecolor=green,linewidth=1pt](10,0){0.5}{0}{180}
 \psarc[linecolor=green,linewidth=1pt](6,0){0.5}{0}{180}
 \psarc[linecolor=blue,linewidth=1pt](6.5,0){0.8}{0}{180}
 \psarc[linecolor=blue,linewidth=1pt](8.5,0){1.2}{0}{180}
\psarc[linecolor=blue,linewidth=1pt](10.5,0){0.8}{0}{180}
 \psarc[linecolor=red,linewidth=1pt](8.5,0){2.8}{0}{180}
  \psarc[linecolor=blue,linewidth=1pt](7.9,0){3.4}{0}{180}
   \psarc[linecolor=blue,linewidth=1pt](5.1,0){0.6}{0}{180}
\psdot(5.7,0)  \rput(5.7,-0.3){$a_1$}
\psdot(7.3,0) \rput(7.3,-0.3){$a_2$}
 \psdot(9.7,0)  \rput(9.7,-0.3){$b_1$}
 \psdot(11.3,0)  \rput(11.3,-0.3){$b_2$}
  \psdot(4.5,0)  \rput(4.5,-0.3){$b_3$}
\psarc[linecolor=red,linewidth=1pt](9.3,0){2}{0}{180}
\end{pspicture}
\caption{The new fundamental domain is the blue one. The ideal triangulation is obtained by the red geodesic. The old fundamental domain is shown in green.}
\end{center}\label{ideal-triang}
\end{figure}

Our new fundamental domain and its  ideal triangulation will be
determined by all the right (or all the left) points of the images
of the bottleneck curve in $\mathbb H$  laying on the absolute. In
figure \ref{pic:bottleneck}, these are $a_1,a_2,b_1,b_2$ and all
their images under the action of the Fuchsian group $\Delta_{g,s}$.
The fundamental domain and the ideal triangulation associated in
this way to our torus with one hole are shown in figure
\ref{ideal-triang}.


\begin{figure}[h]
{\psset{unit=0.5}
\begin{pspicture}(-9,-7)(7,5)
\psellipse[linecolor=cyan](0,-4.7)(1.15,0.5)
\psframe[linecolor=white, fillstyle=solid,
fillcolor=white](-1.2,-4.7)(1.2,-4)
\psellipse[linestyle=dashed,
linecolor=cyan, linewidth=0.5pt](0,-4.7)(1.15,0.5)
\psbezier[linewidth=0.5pt](-2.5,2)(-.5,0.35)(.5,0.35)(2.5,2)
\psbezier[linewidth=0.5pt](-1.5,1.35)(-.5,2.3)(.5,2.3)(1.5,1.35)
\psbezier[linewidth=0.5pt](-3.7,-1.4)(-2.6,-2.)(0.4,-4)(-2.2,-7)
\psbezier[linewidth=0.5pt](-3.7,-1.4)(-6.5,0)(-5.5,5)(0,5)
\psbezier[linewidth=0.5pt](3.7,-1.4)(2.6,-2.)(-0.4,-4)(2.2,-7)
\psbezier[linewidth=0.5pt](3.7,-1.4)(6.5,0)(5.5,5)(0,5)
\psbezier[linecolor=green,
linewidth=0.5pt](0,0.75)(0.5,0.75)(1.1,0.4)(1.1,-0.7)
\psbezier[linecolor=green,
linewidth=0.5pt](1.1,-0.7)(1.1,-2.3)(0.2,-3)(-1.3,-3.7)
\psbezier[linecolor=green,
linewidth=0.5pt](-1.2,-3.65)(-1.5,-3.8)(-1.25,-3.9)(-1.15,-3.75)
\psbezier[linecolor=green, linestyle=dashed,
linewidth=0.5pt](-1.25,-3.8)(-0.5,-3.4)(1,-3.8)(1.2,-4.2)
\psbezier[linecolor=green,
linewidth=0.5pt](1.35,-3.85)(0.75,-4.9)(-1,-4.7)(-1.2,-4.3)
\psbezier[linecolor=green,
linewidth=0.5pt](1.2,-4.3)(0.9,-4.8)(0,-4.9)(-0.7,-4.75)
\psbezier[linecolor=green, linestyle=dashed,
linewidth=0.5pt](-1.2,-4.3)(-0.9,-3.8)(0,-3.9)(0.7,-4.05)
\psbezier[linecolor=green, linestyle=dashed,
linewidth=0.5pt](0,0.75)(-0.5,0.75)(-1.1,0.4)(-1.1,-0.7)
\psbezier[linecolor=green, linestyle=dashed,
linewidth=0.5pt](-1.1,-0.7)(-1.1,-2.3)(-0.2,-3)(1.1,-3.6)
\psbezier[linecolor=green, linestyle=dashed,
linewidth=0.5pt](1.35,-3.85)(1.45,-3.75)(1.2,-3.65)(1,-3.55)
\end{pspicture}
}
\caption{An example of a bounding geodesic asymptotically approaching the
bottleneck curve.}
\label{pic:leo}
\end{figure}
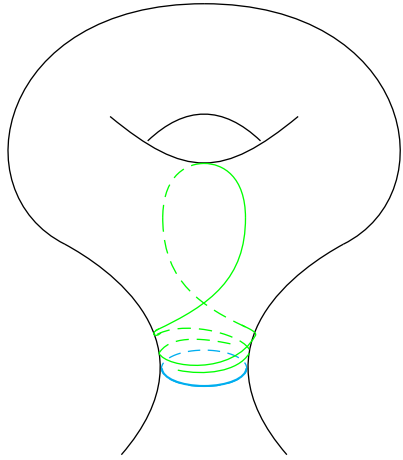

To convince oneself that this new fundamental domain actually
corresponds to the Riemann surface in which all the infinite
hyperboloids attached to the holes have been cut off, consider the
geodesic between $a_2$ and $b_1$. This becomes an infinite (because
$a_2, b_1$ lie on the absolute) geodesic on the Riemann surface that
does not intersect the bottleneck curve (because the latter is not
contained in the new fundamental domain) and therefore winds
infinitely asymptotically approaching the bottleneck curve as in
figure \ref{pic:leo}.


The new fundamental domain is then the finite collection of ideal
triangles (for $F_{g,s}$ we have exactly $4g-4+2 s$ ideal
triangles), which become, under the Fuchsian group action, the
one-to-one covering of an {\em open} Riemann surface with all infinite
hyperboloids removed (together with their bounding bottleneck curves,
to be precise).

By the action of $\mathbb P SL(2,\mathbb R)$, we can move $a_1, a_2$
and $b_2$ to $0,-1$ and $\infty$ respectively. So we can always
label three of the points on the right side of each copy of the
bottleneck curve by $-1,0,\infty$. The fourth will be $\exp(z_0)$
for some real number $z_0$. The element $\gamma_B$ then maps
$a_2\to a_1$, $b_1\to b_2$, $b_2\to b_3$, that is $0\to -1,\quad
\exp(z_0) \to \infty,\quad \infty\to -1-\exp(-z_1),$ as in example
\ref{ex:gammab}.

Let us use the Poincar\'e disk model. We start from the initial
triangle of vertices $0,e^{z_0},\infty$. The transformation
$X_{z_0}$ maps $\infty$ to $0$, $e^{z_0}$ to $-1$ and $0$ to
$\infty$ (see figure \ref{pic:disk-tri}). That is, the initial triangle is mapped in the reference
one. Now we want to reach the final triangle by the same type of
transformation. To this aim, we need to rotate the reference
triangle in clock wise direction. This is represented by the $L$
matrix, or ''left" matrix. In fact $L$ maps $\infty$ to $0$, $0$ to $-1$ and $-1$ to
$\infty$. Now apply $X_{z_1}$, i.e. $0$ is mapped to $\infty$, $-1$
to $e^{z_1}$ and $\infty$ to $0$. Finally we need to apply a right
rotation $R$ to get back to the initial configuration of the reference
triangle.

\begin{figure}[h]
\begin{center}
\begin{pspicture}(4,-4.4)(12,2.4)
  \psarc[linecolor=black,linewidth=2pt](8,-1){3}{0}{360}
  \psarc[linecolor=red,linewidth=2pt](12,-3){3}{112}{195}
  \psarc[linecolor=red,linewidth=2pt](11.2,1.8){2}{181}{260}
\psarc[linecolor=red,linewidth=2pt](16.7,-1.2){8}{160}{198}
\psarc[linecolor=blue,linewidth=2pt](8.27,2.2){1}{194}{336}
\psarc[linecolor=blue,linewidth=2pt](1.6,4.2){8}{300}{342.5}
\psarc[linecolor=blue,linewidth=2pt](3.14,0.84){4.3}{304}{375}
  \psarc[linecolor=black,linewidth=2pt](6.8,-5.07){2.6}{30}{118}
  \rput(9.2,-4){$0$} \rput(11.2,0){$e^{z_0}$}  \rput(9.4,2){$\infty$}  \rput(7.2,2.2){$-1-e^{-z_1}$}
      \rput(5.2,-2.8){$-1$}
      \rput(7.8,-1.8){reference} \rput(7.8,-2.14){triangle}
      \rput(9.5,0){initial}      \rput(9.5,-0.32){triangle}
      \rput(8,1){final}      \rput(8.1,0.7){triangle}
\end{pspicture}
\caption{} \label{pic:disk-tri}
\end{center}
\end{figure}

\begin{figure}[h]
\begin{center}
\begin{pspicture}(1,-1)(16,4)
 \psline{->}(2,0)(12,0)
 \psarc[linecolor=red,linewidth=1pt](8,0){1}{0}{180}
 \psarc[linecolor=blue,linewidth=1pt](8,0){3}{0}{180}
 \psarc[linecolor=blue,linewidth=1pt](7,0){4}{0}{180}
  \psarc[linecolor=blue,linewidth=1pt](4,0){1}{0}{180}
 \psarc(6,0){1}{0}{180}
 \psarc[linecolor=red,linewidth=1pt](10,0){1}{0}{180}
 \rput(5,-0.5){$-1$}
 \rput(7,-0.5){$0$}
 \rput(9,-0.5){$e^{z_0}$}
 \rput(11,-0.5){$\infty$}
 \rput(3,-0.5){$-1-e^{-z_1}$}
 \psarc[linecolor=red,linewidth=1pt](9,0){2}{0}{180}
 \psline[linecolor=green,linewidth=1pt](9,1.6)(8,2.4)
  \psline[linecolor=green,linewidth=1pt](8.8,1.4)(7.8,2.2)
   \psline[linecolor=green,linewidth=1pt](8,2.4)(8,3.4)
    \psline[linecolor=green,linewidth=1pt](7.7,2.4)(7.7,3.4)
     \psline[linecolor=green,linewidth=1pt](7.7,2.4)(6,0.8)
          \psline[linecolor=green,linewidth=1pt](7.9,2.2)(6.2,0.6)
\end{pspicture}
\caption{} \label{pic:ideal1}
\end{center}
\end{figure}

Opening up the same picture on $\mathbb H$ we get again that the red
triangle is mapped into the blue one through the reference one in
the middle (see figure \ref{pic:ideal1}).

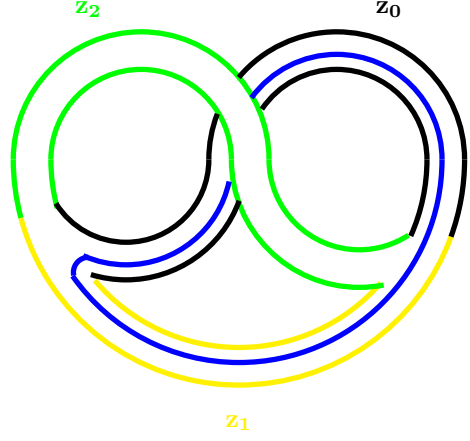
\begin{figure}[h]
\begin{center}
\begin{pspicture}(4,-3.5)(12,2)
 \psarc[linecolor=yellow,linewidth=2pt](8,0){3}{195}{340}
 \psarc[linecolor=yellow,linewidth=2pt](8,0){2.5}{220}{319}
 \psarc[linecolor=black,linewidth=2pt](8,0){2.5}{336}{360}
 \psarc[linecolor=green,linewidth=2pt](8,0){2.5}{180}{195}
  \rput(8,-3.5){\textcolor{yellow}{$\bf{z_1}$}}
 \psarc[linecolor=green,linewidth=2pt](6.7,0){1.7}{0}{180}
  \psarc[linecolor=green,linewidth=2pt](8,0){3}{180}{195}
  \psarc[linecolor=green,linewidth=2pt](6.7,0){1.2}{0}{180}  \psarc[linecolor=green,linewidth=2pt](9.6,0){1.7}{180}{281}
  \psarc[linecolor=green,linewidth=2pt](9.6,0){1.2}{180}{303}
  \rput(6,2){\textcolor{green}{$\bf{z_2}$}}
   \psarc[linecolor=black,linewidth=2pt](8,0){3}{340}{360}
      \psarc[linecolor=black,linewidth=2pt](9.3,0){1.7}{0}{140}
  \psarc[linecolor=black,linewidth=2pt](9.3,0){1.2}{0}{146}
     \psarc[linecolor=black,linewidth=2pt](9.3,0){1.7}{159}{180}
\psarc[linecolor=black,linewidth=2pt](6.5,0){1.6}{253}{340}
  \psarc[linecolor=black,linewidth=2pt](6.5,0){1.1}{212}{360}
 \rput(10,2){\textcolor{black}{$\bf{z_0}$}}
    \psarc[linecolor=blue,linewidth=2pt](8,0){2.7}{215}{360}
         \psarc[linecolor=blue,linewidth=2pt](9.3,0){1.4}{0}{144}
   \psarc[linecolor=blue,linewidth=2pt](6.5,0){1.4}{246}{348}
      \psarc[linecolor=blue,linewidth=2pt](6.0,-1.5){0.2}{90}{192}
\end{pspicture}
\caption{Prezzle with the geodesic corresponding to $\gamma_B$.}\label{pic:gammab}
\end{center}
\end{figure}

We can picture the same decomposition in the fat graph by assigning
a label $z_0,z_1,z_2$ to each edge of the fat graph and choosing an
orientation a priori. For example in figure \ref{pic:gammab} the
decomposition of $\gamma_B$ can be guessed by observing that
starting from the $z_0$ edge of the fatgraph, the curve goes left
into the $z_1$ edge and then right back into the $z_0$ edge.

\begin{figure}[h]
\begin{center}
\begin{pspicture}(4,-3.5)(12,2)
 \psarc[linecolor=yellow,linewidth=2pt](8,0){3}{195}{340}
 \psarc[linecolor=yellow,linewidth=2pt](8,0){2.5}{220}{319}
 \psarc[linecolor=black,linewidth=2pt](8,0){2.5}{336}{360}
 \psarc[linecolor=green,linewidth=2pt](8,0){2.5}{180}{195}
  \rput(8,-3.5){\textcolor{yellow}{$\bf{z_1}$}}
 \psarc[linecolor=green,linewidth=2pt](6.7,0){1.7}{0}{180}
  \psarc[linecolor=green,linewidth=2pt](6.7,0){1.2}{0}{180}  \psarc[linecolor=green,linewidth=2pt](9.6,0){1.7}{180}{281}
  \psarc[linecolor=green,linewidth=2pt](9.6,0){1.2}{180}{303}
  \rput(6,2){\textcolor{green}{$\bf{z_2}$}}
    \psarc[linecolor=green,linewidth=2pt](8,0){3}{180}{195}
       \psarc[linecolor=black,linewidth=2pt](8,0){3}{340}{360}
      \psarc[linecolor=black,linewidth=2pt](9.3,0){1.7}{0}{140}
  \psarc[linecolor=black,linewidth=2pt](9.3,0){1.2}{0}{146}
     \psarc[linecolor=black,linewidth=2pt](9.3,0){1.7}{159}{180}
\psarc[linecolor=black,linewidth=2pt](6.5,0){1.6}{253}{340}
  \psarc[linecolor=black,linewidth=2pt](6.5,0){1.1}{212}{360}
 \rput(10,2){\textcolor{black}{$\bf{z_0}$}}
    \psarc[linecolor=red,linewidth=2pt](9.3,-0.1){1.5}{-56}{140}
   \psarc[linecolor=red,linewidth=2pt](6.7,0){1.5}{0}{180}
   \psline[linecolor=red,linewidth=2pt](8.2,0)(8.25,-0.5)
    \psarc[linecolor=red,linewidth=2pt](9.5,-0.2){1.3}{193}{300}
    \psarc[linecolor=red,linewidth=2pt](6.62,0.1){1.43}{183}{336}
\end{pspicture}
\caption{Prezzle with the geodesic corresponding to $\gamma_A$.}\label{pic:gammaa}
\end{center}
\end{figure}
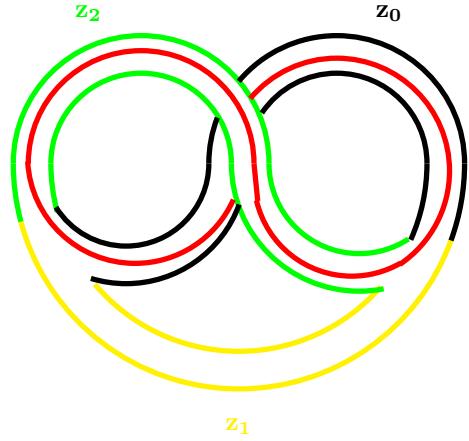

In a similar manner we can obtain also $\gamma_A= L X_{z_2} R X_{z_0}$, see figure  \ref{pic:gammaa}.

Observe that ${\rm dim}({\mathcal T}_1^1)=2$, therefore there must be a relation between $z_0$, $z_1$ and $z_0$.  This relation is
$$
z_0+z_1+z_2=l_P
$$
where $l_P$ is the length of the bottleneck curve $P$. In fact
$$
{\rm Tr} (L X_{z_2} L X_{z_1} L X_{z_0} L X_{z_2} L X_{z_1} L X_{z_0} ) = 2 {\rm cosh}(z_0+z_1+z_2),
$$
and drawing it in the fat--graph we see that this loop unwinds
itself so that it is homotopic to the hole.

 In the general case one can prove the following:

\begin{theorem}\label{param}
Let $F_g^s=\mathbb H\slash\Delta_{g,s}$ and consider ANY
associated fat--graph $\Gamma_g^s$. One can construct all elements
of $\Delta_{g,s}$  by fixing a special starting edge in the fat--graph
and considering all closed paths starting and finishing on the fixed
edge. Then  $\gamma\in \Delta_{g,s}=P_{z_1 z_2 \dots z_n}$ where
$P_{z_1,\dots z_n}$ is a $2\times2$ matrix of the form
$$
P_{z_1 z_2 \dots z_n}= R^{k_n} X_{z_n} R^{k_{n-1}}  \dots X_{z_2} R^{k_1} X_{z_1}
$$
where $k_1, k_2,\dots,k_n =1,2$. Moreover
$$
2 {\rm cosh}(l_\gamma)={\rm Tr}(P_{z_1\dots z_n}),
$$
and $z_1,\dots,z_n$ are coordinates in the Teichm\"uller space extended by $\mathbb R^s$
$$
\mathcal T^s_g\otimes\mathbb R^s.
$$
These coordinates satisfy relations of the form
$$
{\rm Tr}(P_j) = \sum_{i\in I_j} z_i,\qquad j=1,\dots s,
$$
where $I_j$ is the subset of indices labelling all edges in the face
containing the $j$-th hole.
\end{theorem}

\begin{remark}
The number of edges in a fat--graph is $6 g-6 +3 s$ which is the
dimension of $\mathcal T^s_g\otimes\mathbb R^s$. For the $i$-th edge
we have a label $z_i$, so that we have $6 g-6 +3 s$ edges satisfying
$s$ relations. These are enough to determine $\Delta_{g,s}$
uniquely. In fact, if we know for example that the closed path
corresponding to $\gamma$ passes through the edges $x$ and $y$, we
also know $k_x$ and $k_y$ because each vertex is trivalent, so to go
from $x$ to $y$ there only one possible turn.
\end{remark}

\begin{example}
We present here the fat graph corresponding to a  Riemann surface of genus $2$
with two holes.  The geodesic drawn in the picture is uniquely determined by the vertical edges $i, j$
it passes trough. To count the number of faces one is to start on one of the boundaries
and follow it. The number of faces is given by the number of distinct closed
loops we can construct in this way.
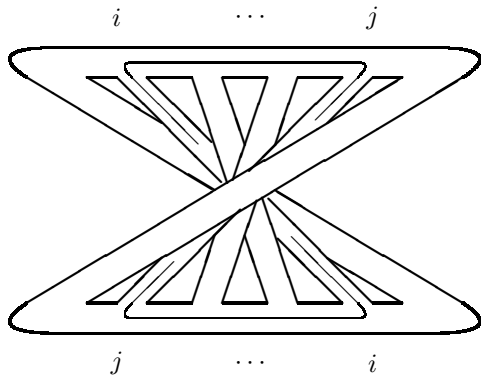
\begin{figure}[h] 
\begin{center}
\setlength{\unitlength}{.4mm}%
\begin{picture}(150,120)(15,-15)
\thicklines
\put(0,95){\line(1, 0){125}}
\qbezier(0,95)(-25,95)(-10,85)
\qbezier(125,95)(150,95)(135,85)
\put(10,85){\line(1, 0){10}} \put(30,85){\line(1, 0){15}}
\put(55,85){\line(1, 0){15}} \put(80,85){\line(1, 0){15}}
\put(105,85){\line(1, 0){10}} \put(0,0){\line(1, 0){125}}
\qbezier(0,0)(-25,0)(-10,10)
\qbezier(125,0)(150,0)(135,10)
\put(10,10){\line(1, 0){10}} \put(30,10){\line(1, 0){15}}
\put(55,10){\line(1, 0){15}} \put(80,10){\line(1, 0){15}}
\put(105,10){\line(1, 0){10}} \put(-10,10){\line(5, 3){125}}
\put(10,10){\line(5, 3){125}} \put(20,10){\line(1, 1){15}}
\put(30,10){\line(1, 1){31}} \put(95,85){\line(-1, -1){31}}
\put(105,85){\line(-1, -1){15}} \put(45,10){\line(1, 3){7.5}}
\put(55,10){\line(1, 3){11.4}} \put(70,85){\line(-1, -3){11.4}}
\put(80,85){\line(-1, -3){7.5}} \put(-10,85){\line(5, -3){62.5}}
\put(10,85){\line(5, -3){25}} \put(115,10){\line(-5, 3){25}}
\put(135,10){\line(-5, 3){62.5}} \put(20,85){\line(1, -1){34.6}}
\put(30,85){\line(1, -1){21.5}} \put(95,10){\line(-1, 1){21.5}}
\put(105,10){\line(-1, 1){34.6}} \put(45,85){\line(1, -3){11.7}}
\put(55,85){\line(1, -3){7.5}} \put(70,10){\line(-1, 3){7.5}}
\put(80,10){\line(-1, 3){11.7}}
\thinlines
\qbezier(25,90)(20,90)(25,85)
\qbezier(100,90)(105,90)(100,85)
\qbezier(25,5)(20,5)(25,10)
\qbezier(100,5)(105,5)(100,10)
\put(25,90){\line(1,0){75}}
\put(25,5){\line(1,0){75}}
\put(25,85){\line(1,-1){22}}
\put(100,10){\line(-1,1){22}}
\put(100,85){\line(-1,-1){14}}
\put(25,10){\line(1,1){14}}
\put(20,105){\makebox(0,0){$i$}}
\put(65,105){\makebox(0,0){$\cdots$}}
\put(105,105){\makebox(0,0){$j$}}
\put(20,-10){\makebox(0,0){$j$}}
\put(65,-10){\makebox(0,0){$\cdots$}}
\put(105,-10){\makebox(0,0){$i$}}
\end{picture}
\label{octopus}\caption{Fat graph corresponding to a  Riemann surface of genus $2$
with two holes.}
\end{center}
\end{figure}

\end{example}

\section{Modular group action}

Let us draw the attention now to the word ``ANY" in the statement of
Theorem \ref{param}. This corresponds to the fact that our splitting
of the crossings in the fat--graph was arbitrary. Choosing a
different splitting corresponds to perform a transformation that
preserves all geodesics and all traces, i.e. to the action of the
modular group $\mathbb{P}SL(2,\mathbb R)$.

On the fat-graph changing the splitting of the crossing means that the edge which was
on top, now goes on the bottom and viceversa, as in figure \ref{pic:flip}

\begin{figure}[h]
\begin{center}
\begin{pspicture}(1.2,-3.5)(4.2,2)
  \psarc[linecolor=yellow,linewidth=2pt](8,0){3}{195}{340}
 \psarc[linecolor=yellow,linewidth=2pt](8,0){2.5}{220}{319}
 \psarc[linecolor=black,linewidth=2pt](8,0){2.5}{336}{360}
 \psarc[linecolor=green,linewidth=2pt](8,0){2.5}{180}{195}
  \rput(8,-3.5){\textcolor{yellow}{$\bf{z_0}$}}
    \psarc[linecolor=green,linewidth=2pt](8,0){3}{180}{195}
 \psarc[linecolor=green,linewidth=2pt](6.7,0){1.7}{40}{180}
  \psarc[linecolor=green,linewidth=2pt](6.7,0){1.7}{0}{20}
  \psarc[linecolor=green,linewidth=2pt](6.7,0){1.2}{34}{180}
   \psarc[linecolor=green,linewidth=2pt](9.6,0){1.7}{198}{281}
  \psarc[linecolor=green,linewidth=2pt](9.6,0){1.2}{180}{303}
  \rput(6,2){\textcolor{green}{$\bf{z_2}$}}
  \psarc[linecolor=black,linewidth=2pt](8,0){3}{340}{360}
      \psarc[linecolor=black,linewidth=2pt](9.3,0){1.7}{0}{180}
  \psarc[linecolor=black,linewidth=2pt](9.3,0){1.2}{0}{180}
\psarc[linecolor=black,linewidth=2pt](6.5,0){1.6}{253}{360}
  \psarc[linecolor=black,linewidth=2pt](6.5,0){1.1}{212}{360}
 \rput(10,2){\textcolor{black}{$\bf{z_1}$}}
  \psline{->}(4.2,-1)(4.7,-1)
 \end{pspicture}
\begin{pspicture}(14.8,-3.5)(11.4,2)
 \psarc[linecolor=yellow,linewidth=2pt](8,0){3}{195}{340}
 \psarc[linecolor=yellow,linewidth=2pt](8,0){2.5}{220}{319}
 \psarc[linecolor=black,linewidth=2pt](8,0){2.5}{336}{360}
 \psarc[linecolor=green,linewidth=2pt](8,0){2.5}{180}{195}
  \rput(8,-3.5){\textcolor{yellow}{$\bf{\tilde z_0}$}}
   \psarc[linecolor=green,linewidth=2pt](8,0){3}{180}{195}
 \psarc[linecolor=green,linewidth=2pt](6.7,0){1.7}{0}{180}
  \psarc[linecolor=green,linewidth=2pt](6.7,0){1.2}{0}{180}  \psarc[linecolor=green,linewidth=2pt](9.6,0){1.7}{180}{281}
  \psarc[linecolor=green,linewidth=2pt](9.6,0){1.2}{180}{303}
  \rput(6,2){\textcolor{green}{$\bf{\tilde z_2}$}}
   \psarc[linecolor=black,linewidth=2pt](8,0){3}{340}{360}
      \psarc[linecolor=black,linewidth=2pt](9.3,0){1.7}{0}{140}
  \psarc[linecolor=black,linewidth=2pt](9.3,0){1.2}{0}{146}
     \psarc[linecolor=black,linewidth=2pt](9.3,0){1.7}{159}{180}
\psarc[linecolor=black,linewidth=2pt](6.5,0){1.6}{253}{340}
  \psarc[linecolor=black,linewidth=2pt](6.5,0){1.1}{212}{360}
 \rput(10,2){\textcolor{black}{$\bf{\tilde z_1}$}}
 \end{pspicture}
 \caption{Flip}\label{pic:flip}
 \end{center}
\end{figure}

We can visualize this edge flip by representing the prezzle itself
as a graph. In fact as we saw in figure \ref{pic:torus-to-pr}, the
fat graph is equivalent to a choice of the way of resolving the
crossing of the  bouquet of cycles associated to our Riemann
surface. Schematically we can visualize such resolving of the
crossing as in figure \ref{pic:crossing}.

\begin{figure}[h]
\begin{pspicture}(7,-1.5)(11,2)
 \pscurve[linecolor=red,linewidth=1.5pt,showpoints=false](5.2,.4)(5.1,1.2)(6,1.3)(6.3,0.3)(7,0)
 \pscircle[fillstyle=solid,fillcolor=white,linecolor=white](6.2,0.75){0.15}
 \psellipse[linecolor=blue,linewidth=1.5pt](6,0)(1,0.8)
 \psline{->}(8.5,0)(9,0)
  \end{pspicture}
 \begin{pspicture}(6,-1.5)(9,2)
  \pscurve[linecolor=red,linewidth=1.5pt,showpoints=false](5.2,.4)(5.1,1.2)(6,1.3)(6.3,0.3)(7,0)
   \pscircle[fillstyle=solid,fillcolor=white,linecolor=white](6.0,1){0.15}
 \pscircle[fillstyle=solid,fillcolor=white,linecolor=white](6.2,0.75){0.5}
  \pscircle[fillstyle=solid,fillcolor=white,linecolor=white](5.5,1.5){0.7}
 \psellipse[linecolor=blue,linewidth=1.5pt](6,0)(1,0.8)
 \pscircle[fillstyle=solid,fillcolor=white,linecolor=white](5.7,-1){1}
  \psline{->}(8,0)(8.5,0)
   \end{pspicture}
  \begin{pspicture}(-1.,-2.4)(-0.2,0)
   \psline[linecolor=red,linewidth=1.5pt](-1,-3)(0,-2)\rput(-1.5,-3.3){$z_1$}
 \psline[linecolor=blue,linewidth=1.5pt](0,-2)(1,-3) \rput(1.5,-3.3){$z_2$}
  \psline[linecolor=blue,linewidth=1.5pt](0,-2)(0,0) \rput(0.5,-1){$z_0$}
 \psline[linecolor=blue,linewidth=1.5pt](-1,1)(0,0)
 \psline [linecolor=red,linewidth=1.5pt](0,0)(1,1) \rput(-1.5,1.3){$z_2 $} \rput(1.5,1.3){$z_1$}
\end{pspicture}
\caption{}\label{pic:crossing}
 \end{figure}

A different resolving of the crossing leads to another prezzle as in
figure \ref{pic:flip} and to another graph as in figure
\ref{pic:crossing1}.

\begin{figure}[h]
\begin{pspicture}(8,-1.5)(11,2)
  \psellipse[linecolor=red,linewidth=1.5pt,showpoints=false](6.2,0.5)(1,0.8)
   \pscircle[fillstyle=solid,fillcolor=white,linecolor=white](6.2,0.75){1}
     \pscircle[fillstyle=solid,fillcolor=white,linecolor=white](7.2,0.1){.5}
          \pscircle[fillstyle=solid,fillcolor=white,linecolor=white](5.2,0.6){.15}
\psarc[linecolor=red,linewidth=1.5pt,showpoints=false](6.75,0.75){.7}{310}{180}
 \pscircle[fillstyle=solid,fillcolor=white,linecolor=white](6.2,0.75){0.15}
 \psellipse[linecolor=blue,linewidth=1.5pt](6,0)(1,0.8)
 \psline{->}(8.5,0)(9,0)
  \end{pspicture}
 \begin{pspicture}(6,-1.5)(9,2)
  \psellipse[linecolor=red,linewidth=1.5pt,showpoints=false](6.2,0.5)(1,0.8)
   \pscircle[fillstyle=solid,fillcolor=white,linecolor=white](6.2,0.75){1}
     \pscircle[fillstyle=solid,fillcolor=white,linecolor=white](7.2,0.1){.5}
          \pscircle[fillstyle=solid,fillcolor=white,linecolor=white](5.2,0.6){.15}
\psarc[linecolor=red,linewidth=1.5pt,showpoints=false](6.75,0.75){.7}{310}{180}
 \pscircle[fillstyle=solid,fillcolor=white,linecolor=white](6.2,0.75){0.15}
 \psellipse[linecolor=blue,linewidth=1.5pt](6,0)(1,0.8)
    \pscircle[fillstyle=solid,fillcolor=white,linecolor=white](7.3,0.75){.8}
       \pscircle[fillstyle=solid,fillcolor=white,linecolor=white](6.2,-0.6){1}
  \psline{->}(7.6,0)(8.1,0)
   \end{pspicture}
  \begin{pspicture}(11,-2.5)(9,0)
  \psline[linecolor=red,linewidth=1.5pt](8,-2)(9,-1)
  \psline[linecolor=blue,linewidth=1.5pt](9,-1)(8,0)
  \rput(8,0.5){$\tilde z_2$} \rput(8,-2.5){$\tilde z_1$}
 \psline[linecolor=blue,linewidth=1.5pt](9,-1)(11,-1) \rput(10,-1.5){$\tilde z_0$}
 \psline[linecolor=blue,linewidth=1.5pt](12,-2)(11,-1)
 \psline[linecolor=red,linewidth=1.5pt](11,-1)(12,0) \rput(12.4,0.5){$\tilde z_1$} \rput(12.4,-2.5){$\tilde z_2$}
 \end{pspicture}
\caption{}\label{pic:crossing1}
 \end{figure}

Observe that the torus example is special in the sense that each
couple of edges intersect twice. In general we'll have $5$ different
labels $z_0,z_1,z_2,z_3,z_4$. The labels are by convention fixed
like in figure \ref{pic:hyp1}.

To find the way in which the labels are transformed, it is useful to
realize that a different splitting of the crossing corresponds to a
different ideal triangulation of the hyperbolic plane in triangles
as described in figure \ref{pic:hyp1}.

\begin{figure}[h]
\begin{center}
\begin{pspicture}(2,-0.6)(16,3.4)
 \psline{->}(4,0)(12,0)
 \psarc(8,0){1}{0}{180}
 \psarc(8,0){3}{0}{180}
 \psarc(6,0){1}{0}{180}
 \psarc(10,0){1}{0}{180}
 \rput(5,-0.5){$x_1$}
 \rput(7,-0.5){$x_2$}
 \rput(9,-0.5){$x_3$}
 \rput(11,-0.5){$x_4$}
 \psarc[linecolor=red,linewidth=1pt](9,0){2}{0}{180}
  \psline[linecolor=green,linewidth=1pt](8,1.3)(10,0.5)  \rput(10.2,.4){$z_3$}
 \psline[linecolor=green,linewidth=1pt](8,1.3)(7.8,0.5)  \rput(7.8,0.3){$z_2$}
  \psline[linecolor=green,linewidth=1pt](8,1.3)(8,2.4) \rput(8.2,2){$z_0$}
    \psline[linecolor=green,linewidth=1pt](8,2.4)(10,3)  \rput(10.2,3){$z_1$}
      \psline[linecolor=green,linewidth=1pt](8,2.4)(6,.4) \rput(5.7,.3){$z_4$}
\end{pspicture}
\begin{pspicture}(2,-0.6)(16,3.4)
 \psline{->}(4,0)(12,0)
 \psarc(8,0){1}{0}{180}
 \psarc(8,0){3}{0}{180}
 \psarc(6,0){1}{0}{180}
 \psarc(10,0){1}{0}{180}
 \rput(5,-0.5){$x_1$}
 \rput(7,-0.5){$x_2$}
 \rput(9,-0.5){$x_3$}
 \rput(11,-0.5){$x_4$}
 \psarc[linecolor=blue,linewidth=1pt](7,0){2}{0}{180}
      \psline[linecolor=green,linewidth=1pt](7.8,1.5)(9,1.5)  \rput(8.7,1.8){$\tilde z_0$}
\psline[linecolor=green,linewidth=1pt](9,1.5)(10,0.5)  \rput(10.2,.35){$\tilde z_3$}
     \psline[linecolor=green,linewidth=1pt](9,1.5)(9.7,3.2)  \rput(9.9,3.1){$\tilde z_1$}
\psline[linecolor=green,linewidth=1pt](7.8,1.5)(7.6,0.5)  \rput(7.7,0.3){$\tilde z_2$}
\psline[linecolor=green,linewidth=1pt](7.8,1.5)(6,0.5)  \rput(5.8,0.4){$\tilde z_4$}
           \end{pspicture}
\caption{}\label{pic:hyp1}
\end{center}
\end{figure}

\begin{figure}[h]
\begin{center}
\psset{unit=7mm}
\begin{pspicture}(-1.5,-3.5)(13,1)
 \psline(-1,-3)(0,-2)(1,-3) \rput(-1.5,-3.5){$z_2$} \rput(1.5,-3.5){$z_3$}
 \psline(0,-2)(0,0) \rput(0.5,-1){$z_0$}
 \psline(-1,1)(0,0)(1,1) \rput(-1.5,1.5){$z_4$} \rput(1.5,1.5){$z_1$}
 \psline{->}(3.7,-1)(5.7,-1)
 \psline(8,-2)(9,-1)(8,0) \rput(8,0.5){$z_4+\varphi(z_0)$} \rput(8,-2.5){$z_2-\varphi(-z_0)$}
 \psline(9,-1)(11,-1) \rput(10,-1.5){$-z_0 $}
 \psline(12,-2)(11,-1)(12,0) \rput(13,0.5){$z_1-\varphi(-z_0)$} \rput(13,-2.5){$z_3+\varphi(z_0)$}
\end{pspicture}\caption{The flip in the general situation.}\label{gen-flip}
\end{center}
\end{figure}
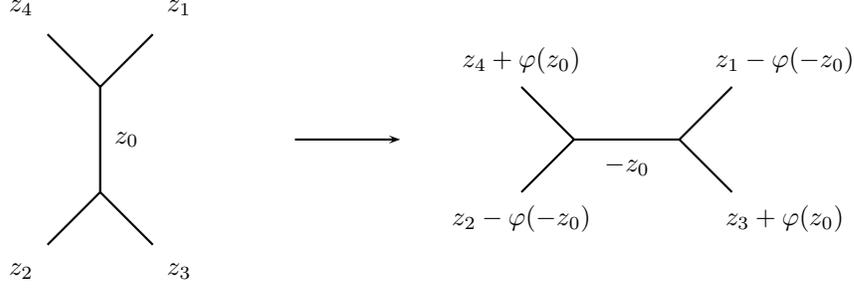

If we denote the vertices $-1,0,e^{z_0},\infty$ by $x_1,\dots,x_4$,
we see that the first splitting gives the {\em four-term relation}\footnote{In the standard
terminology, variables $z_i$ are called the shear coordinates.}
$$
\frac{(x_2-x_1)(x_4-x_3)}{(x_3-x_2)(x_4-x_1)}=e^{-z_0},
$$
while the second gives
$$
\frac{(x_3-x_2)(x_4-x_1)}{(x_4-x_3)(x_2-x_1)}=e^{z_0},
$$
so that the flip corresponds to map  $z_0\to -z_0$. Let us see how
$z_2$ changes. To this aim,  let us pick a fifth point $x_5$, such
that $x_2<x_5<x_3$ and construct the hyperbolic triangle with
vertices $x_2,x_5,x_3$. If we look under the microscope, we see that
the initial configuration for $z_2$ is the same as the $z_0$
configuration in the bottom of figure \ref{pic:hyp1} in which the
point labels are changed as
$$
x_1\rightarrow x_2,\quad x_2\rightarrow x_5,
\quad x_3\rightarrow x_3,\quad x_4\rightarrow x_4.
$$
So by the four-term relation we get
$$
\frac{(x_4-x_3)(x_5-x_2)}{(x_3-x_5)(x_4-x_2)}=e^{z_2}.
$$
In the same way, the final configuration for $\tilde z_2$ is the
same as the $z_0$ configuration in the top figure
 \ref{pic:hyp1},  again with points changed as
 $$
x_1\rightarrow x_1,\quad x_2\rightarrow x_2,
\quad x_3\rightarrow x_5,\quad x_4\rightarrow x_3.
$$
So:
$$
\frac{(x_5-x_2)(x_3-x_1)}{(x_2-x_1)(x_3-x_5)}=e^{-\tilde z_2}.
$$
It is a straightforward computation to show that then $\tilde z_2=
z_2+\varphi(z_0)$ where $\varphi(z_0)=\log(1+e^{z_0})$. One can
prove analogous relations for $\tilde z_1, \tilde z_3$ and $\tilde
z_4$. These are represented in figure \ref{gen-flip}.

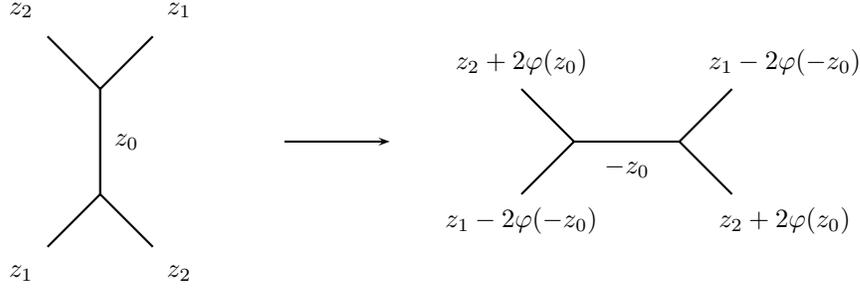
\begin{figure}[h]
\begin{center}
\psset{unit=7mm}
\begin{pspicture}(-1.5,-3.5)(13,1)
 \psline(-1,-3)(0,-2)(1,-3) \rput(-1.5,-3.5){$z_1$} \rput(1.5,-3.5){$z_2$}
 \psline(0,-2)(0,0) \rput(0.5,-1){$z_0 $}
 \psline(-1,1)(0,0)(1,1) \rput(-1.5,1.5){$z_2$} \rput(1.5,1.5){$z_1$}
 \psline{->}(3.5,-1)(5.5,-1)
 \psline(8,-2)(9,-1)(8,0) \rput(8,0.5){$z_2+2\varphi(z_0)$} \rput(8,-2.5){$z_1-2\varphi(-z_0)$}
 \psline(9,-1)(11,-1) \rput(10,-1.5){$-z_0 $}
 \psline(12,-2)(11,-1)(12,0) \rput(13,0.5){$z_1-2\varphi(-z_0)$} \rput(13,-2.5){$z_2+2\varphi(z_0)$}
\end{pspicture}\caption{The flip in the torus situation.}\label{tor-flip}
\end{center}
\end{figure}

In the torus situation, every edges meets every other one twice.
This means that we have to sum the contributions. In fact, for
example $z_2=z_4$. Therefore we have a double contribution $2
\varphi(z)$. This is shown in figure \ref{tor-flip}.

Let us see how to produce a mapping class group transformation from
a flip morphism. We call it ``morphism" because it establishes a
relation between two specific graphs. Denote the morphism between
two graphs $\Gamma_1$ and $\Gamma_2$ by $[\Gamma_1,\Gamma_2]$.
Suppose we have a sequence of morphisms between graphs of the same
combinatorial type (i.e. only the labels are different):
$$
[\Gamma_1,\Gamma_2] ,[\Gamma_2,\Gamma_3],\dots,[\Gamma_{n-1},\Gamma_n],
$$
denote by $z_\alpha^{(k)}$ the labels of the $k$-th graph
$\Gamma_k$. We want to show that $z_\alpha^{(n)}$ is obtained from
$z_{\beta}^{(1)}$ by the action of the mapping class group. Again we
do it in our example ${\mathcal T}_1^1$.

Under this map the order of $z_0,z_1,z_2$ must be preserved: since
before the mapping class group transformation, by a anti-clockwise
rotation, we have the order $z_0,z_1,z_2$, and after the mapping
class group transformation, we have the order $-z_0,
z_2+2\varphi(z),z_1-2\varphi(-z)$, we see that the map is
$$
z_0\mapsto -z_0, \qquad z_1\mapsto z_2+ 2 \log(1+e^{z_0}),\qquad
z_2\mapsto z_1- 2 \log(1+e^{-z_0}),
$$
so that $z_0+z_1+z_2$ remains invariant and the orientation is preserved.

Let  $U=\exp(\frac{z_0}{2})$,  $V=\exp(\frac{z_2}{2})$, then
$$
U\mapsto U^{-1},
$$
and
$$
V\mapsto \frac{e^{\frac{z_1}{2}}}{1+e^{-z_0}}=e^{\frac{l_P}{2}} \frac{e^{\frac{-z_2}{2}}
e^{\frac{-z_0}{2}}}{1+e^{-z_0}}
=e^{\frac{l_P}{2}}\frac{V^{-1}}{U+U^{-1}},
$$
so that indeed we have a modular transformation. This is in fact a general statement.

The next theorem is valid in the general situation of SINGLE intersection (not for the torus!):

\begin{theorem}\cite{pen1,pen2} The following special relations between morphisms hold true:
\begin{enumerate}
\item The flips are in involution:
$$
[\Gamma_1,\Gamma_2]  [\Gamma_2,\Gamma_1]={\rm Id},
$$
\item  Four-term relation: flips of disjoint edges commute.
\item Pentagon identity (valid for flips of two edges having a single common vertex):
\end{enumerate}
\label{th:morph}
\end{theorem}

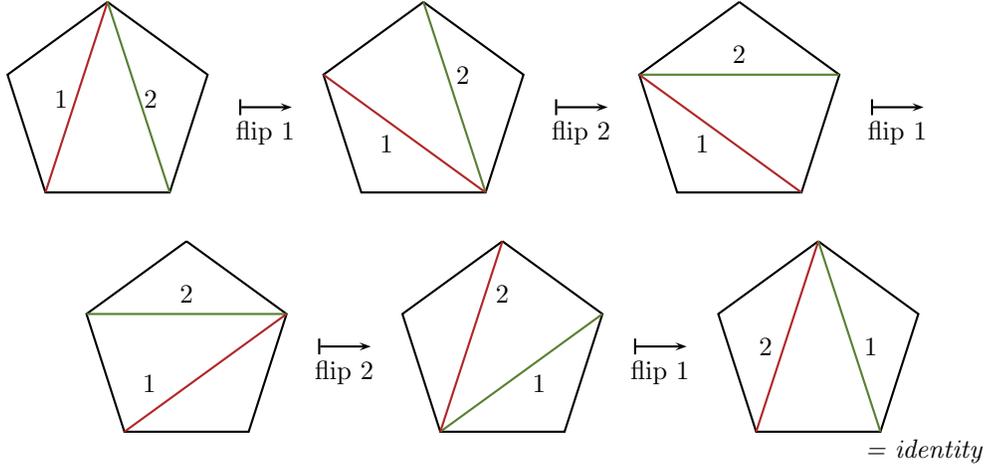
\begin{center}
\newrgbcolor{olivegreen}{0.35 0.5 0.2}
\newrgbcolor{firebrick}{0.7 0.15 0.15}
\psset{unit=0.7}
\begin{pspicture}(-2,-2)(16,2.5)
 \psline(0,2)(-1.90, 0.618)(-1.18, -1.62)(1.18, -1.62)(1.90,0.618)(0,2)
   \psline[linecolor=firebrick](-1.18, -1.62)(0,2)\psline[linecolor=olivegreen](0,2)(1.18, -1.62) \put(-1,0){1} \put(0.7,0){2}
     \psline{|->}(2.5,0)(3.5,0) \rput(3,-0.5){flip 1}
 \rput(6,0){\psline(0,2)(-1.90, 0.618)(-1.18, -1.62)(1.18, -1.62)(1.90,0.618)(0,2)
   \psline[linecolor=olivegreen](0,2)(1.18, -1.62)\psline[linecolor=firebrick](1.18, -1.62)(-1.90,0.618) \rput(-0.7, -0.7){1} \rput(0.75, 0.6){2}
      \psline{|->}(2.5,0)(3.5,0) \rput(3,-0.5){flip 2}
 }
 \rput(12,0){\psline(0,2)(-1.90, 0.618)(-1.18, -1.62)(1.18, -1.62)(1.90,0.618)(0,2)
   \psline[linecolor=olivegreen](1.90, 0.618)(-1.90, 0.618)\psline[linecolor=firebrick](-1.90, 0.618)(1.18, -1.62) \rput(-0.7,-0.7){1}\rput(0,1){2}
      \psline{|->}(2.5,0)(3.5,0) \rput(3,-0.5){flip 1}
 }
\end{pspicture}

\begin{pspicture}(-2,-2)(13,2.5)
 \psline(0,2)(-1.90, 0.618)(-1.18, -1.62)(1.18, -1.62)(1.90,0.618)(0,2)
   \psline[linecolor=olivegreen](-1.90, 0.618)(1.90, 0.618)\psline[linecolor=firebrick](1.90, 0.618)(-1.18, -1.62)  \rput(-0.7,-0.7){1}\rput(0,1){2}
   \psline{|->}(2.5,0)(3.5,0) \rput(3,-0.5){flip 2}
 \rput(6,0){\psline(0,2)(-1.90, 0.618)(-1.18, -1.62)(1.18, -1.62)(1.90,0.618)(0,2)
   \psline[linecolor=olivegreen](1.90, 0.618)(-1.18, -1.62)\psline[linecolor=firebrick](-1.18, -1.62)(0,2)  \rput(0.7,-0.7){1} \rput(0,1){2}
      \psline{|->}(2.5,0)(3.5,0) \rput(3,-0.5){flip 1}
  }
 \rput(12,0){\psline(0,2)(-1.90, 0.618)(-1.18, -1.62)(1.18, -1.62)(1.90,0.618)(0,2)
   \psline[linecolor=firebrick](-1.18, -1.62)(0,2)\psline[linecolor=olivegreen](0,2)(1.18, -1.62)  \rput(1,0){1} \rput(-1,0){2}
  }
 \rput(14,-2){\it = identity}
\end{pspicture}
\end{center}

\section{Poisson brackets}

Since quantization is our ultimate aim, we need to select an
appropriate set of observables. The natural observables of the
theory are the lengths of the closed geodesics. We shall denote a collection of labels
by $P=\{z_{i_1},\dots,z_{i_n}\}$, and identify it with the
corresponding matrix $P=R^{k_n} X_{z_{i_n}} R^{k_{n-1}}\dots
X_{z_{i_2}} R^{k_1} X_{z_{i_1}}$. We'll denote by $G_P$ the length
of the geodesic $\gamma_P$ passing trough the edges
$\{z_{i_1},\dots,z_{i_n}\}$.

We postulate the Poisson brackets between edges as follows: we
establish an ordering, for example anticlockwise, at each edge, so
that we label the edges concurring at the same vertex by $z_1, z_2,
z_3$. Then we establish
\begin{equation}\label{poisson}
\{z_i,z_{i+1}\}=1.
\end{equation}
To find the Poisson bracket between two edges, we need to sum on all common vertices.
In our example of the torus with one hole we get
$$
\{z_0,z_1\}=\{z_1,z_2\}=\{z_2,z_0\}=2.
$$

\begin{theorem}
The only central elements of the given Poisson algebra are ${\rm
Tr}(P_j) = \sum_{i\in I_j} z_i$, $j=1,\dots s,$, where $I_j$ is the
subset of indices labeling all edges in the face containing the
$j$-th hole.
\end{theorem}

\begin{corollary}
The dimension of the symplectic leaves is equal to the dimension of
the Teichm\"uller space, which is $6 g-6+2s$.
\end{corollary}

\subsection{Relations between geodesic functions}
\begin{figure}[h]\begin{center}
\begin{pspicture}(-3,-4)(0.6,1)
  \psellipse[linecolor=blue,linewidth=1pt](-1,-1)(1,2)
 \pscircle[fillstyle=solid,fillcolor=white,linecolor=white](-1.1,-2.9){0.2}
 \rput(-0.2,-3.0){$\textcolor{blue}{G_Q}$}
\psellipse[linecolor=red,linewidth=1pt](-2,-2)(-1,2)\rput(-2.8,-3.8){$\textcolor{red}{G_P}$}
 \pscircle[linecolor=black,linewidth=2pt](-2,-0.2){0.7}\rput(0.5,-1.5){$=$}
\end{pspicture}
\begin{pspicture}(-3,-4)(0.6,1)
  \psellipse[linecolor=blue,linewidth=1pt](-1,-1)(1,2)
 \pscircle[fillstyle=solid,fillcolor=white,linecolor=white](-1.1,-2.9){0.2}
\psellipse[linecolor=red,linewidth=1pt](-2,-2)(-1,2)
 \pscircle[fillstyle=solid,fillcolor=white,linecolor=black,linewidth=2pt](-2,-0.2){0.7}
 \psarc[linewidth=1.5pt](-1,.2){0.75}{166}{240}
 \psarc[linewidth=1.5pt](-2.6,-1.2){0.75}{30}{90}\rput(-0.6,-3.5){$G_{PQ}$}
 \rput(0.45,-1.5){$+$}
\end{pspicture}
\begin{pspicture}(-3,-4)(0,1)
  \psellipse[linecolor=blue,linewidth=1pt](-1,-1)(1,2)
 \pscircle[fillstyle=solid,fillcolor=white,linecolor=white](-1.1,-2.9){0.2}
\psellipse[linecolor=red,linewidth=1pt](-2,-2)(-1,2)
 \pscircle[fillstyle=solid,fillcolor=white,linecolor=black,linewidth=2pt](-2,-0.2){0.7}
\psarc[linewidth=1.5pt](-1.3,-1.22){0.75}{96}{154}
\psarc[linewidth=1.5pt](-2.6,.4){0.9}{268}{0}
\rput(-0.6,-3.5){$G_{PQ^{-1}}$}
\end{pspicture}\caption{Skein relation.}\label{pic:skein}
\end{center}
\end{figure}
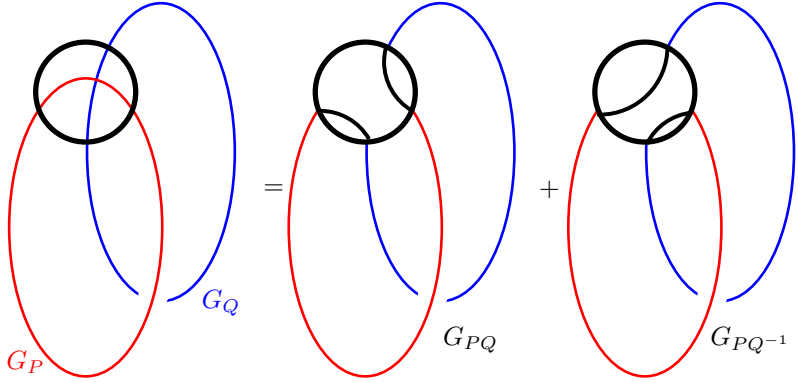

Let us first deal with the case of one intersection.
\begin{itemize}
\item Skein--relation:
$$
G_PG_Q =G_{PQ} + G_{PQ^{-1}}.
$$
This corresponds to the relation
$$
{\rm Tr}(P Q)+{\rm Tr}(P Q^{-1})
={\rm Tr}(P){\rm Tr}(Q),
$$
valid for matrices having determinant equal to one. Graphically it
means that we split the crossing in two possible ways (see figure
\ref{pic:skein}).

\end{itemize}

\begin{itemize}

\item Poisson bracket: it is induced by bracket (\ref{poisson})\footnote{See the calculation in Appendix~B.}
\begin{equation}\label{eq:poisson1}
\{G_P,G_Q\} =\frac{1}{2} G_{PQ} -\frac{1}{2} G_{PQ^{-1}}.
\end{equation}
This is precisely the Goldman bracket on the space of geodesic
functions in $2+1$--gravity. It comes from the assumption
(\ref{poisson}). This is visualized graphically as in figure
\ref{pic:poisson}. To remember this rule we fix an orientation on
$\gamma_P$ and establish to assign a $+$ when we turn left and a
minus when we turn right. Observe that since ${\rm Tr}(P) = {\rm
Tr}(P^{-1})$, this Poisson bracket does not depend on the
orientation of the loops $\gamma_P$ and $\gamma_Q$.
\end{itemize}

\begin{figure}[h]\begin{center}
\begin{pspicture}(-3.3,-4)(0.98,1)
\rput(-3.2,-1.5){$\Bigg\{$}
  \psellipse[linecolor=blue,linewidth=1pt](-1,-1)(1,2)
 \pscircle[fillstyle=solid,fillcolor=white,linecolor=white](-1.1,-2.9){0.2}
  \rput(-0.2,-3.0){$\textcolor{blue}{G_Q}$}
\psellipse[linecolor=red,linewidth=1pt](-2,-2)(-1,2)
\rput(-2.8,-3.8){$\textcolor{red}{G_P}$}
 \pscircle[linecolor=black,linewidth=2pt](-2,-0.2){0.7}\rput(0.5,-1.5){$\Bigg\}=\,\frac{1}{2}$}
\end{pspicture}
\begin{pspicture}(-3,-4)(0.8,1)
  \psellipse[linecolor=blue,linewidth=1pt](-1,-1)(1,2)
 \pscircle[fillstyle=solid,fillcolor=white,linecolor=white](-1.1,-2.9){0.2}
\psellipse[linecolor=red,linewidth=1pt](-2,-2)(-1,2)
 \pscircle[fillstyle=solid,fillcolor=white,linecolor=black,linewidth=2pt](-2,-0.2){0.7}
 \psarc[linewidth=1.5pt](-1,.2){0.75}{166}{240}
 \psarc[linewidth=1.5pt](-2.6,-1.2){0.75}{30}{90}
 \rput(-0.6,-3.5){$G_{PQ}$}
\rput(0.5,-1.5){$-\frac{1}{2}$}
\end{pspicture}
\begin{pspicture}(-3,-4)(0,1)
  \psellipse[linecolor=blue,linewidth=1pt](-1,-1)(1,2)
 \pscircle[fillstyle=solid,fillcolor=white,linecolor=white](-1.1,-2.9){0.2}
\psellipse[linecolor=red,linewidth=1pt](-2,-2)(-1,2)
 \pscircle[fillstyle=solid,fillcolor=white,linecolor=black,linewidth=2pt](-2,-0.2){0.7}
\psarc[linewidth=1.5pt](-1.3,-1.22){0.75}{96}{154}
\psarc[linewidth=1.5pt](-2.6,.4){0.9}{268}{0}
\rput(-0.6,-3.5){$G_{PQ^{-1}}$}
\end{pspicture}
\caption{Poisson bracket.}\label{pic:poisson}
\end{center}
\end{figure}

\begin{example}
Let us show that the Poisson bracket between geodesics coincides
with (\ref{poisson}) in the special case of
$P=\gamma_A=L X_{z_2} R X_{z_0}$ and $Q=\gamma_B = R X_{z_1} L
X_{z_0}$
$$
P=\left(\begin{array}{cc}e^{-\frac{{z_2}}{2}-\frac{{z_0}}{2}}&
-e^{-\frac{{z_2}}{2}+\frac{{z_0}}{2}}\\
 -e^{-\frac{{z_2}}{2}-\frac{{z_0}}{2}}&e^{-\frac{{z_0}}{2}}(e^{\frac{{z_2}}{2}}+e^{-\frac{{z_2}}{2}})\\
\end{array}\right),\qquad
Q=\left(\begin{array}{cc}e^{-\frac{{z_0}}{2}}(e^{-\frac{{z_1}}{2}}+e^{\frac{{z_1}}{2}})&
-e^{\frac{{z_1}}{2}+\frac{{z_0}}{2}}\\
 -e^{\frac{{z_1}}{2}-\frac{{z_0}}{2}}&e^{\frac{{z_1}}{2}+\frac{{z_0}}{2}}\\
\end{array}\right),
$$
\begin{eqnarray}
\nn&&
\left\{{\rm Tr}(P),{\rm Tr}(Q)\right\} =
\frac{\partial {\rm Tr}(P)}{\partial {z_2}} \frac{\partial {\rm Tr}(Q)}{\partial {z_0}}\{{z_2},z_0\}+
\frac{\partial {\rm Tr}(P)}{\partial {z_0}} \frac{\partial {\rm Tr}(Q)}{\partial {z_1}}\{{z_0},{z_1}\}+\nn\\&&
\qquad\qquad\qquad\quad+ \frac{\partial {\rm Tr}(P)}{\partial {z_2}}
\frac{\partial {\rm Tr}(Q)}{\partial {z_1}}\{{z_2},{z_1}\}=\nn\\
&&\nn
=\frac{1}{2} e^{-\frac{z_2}{2}-\frac{{z_1}}{2}-{z_0}}
\left(1+e^{{z_1}}-e^{{z_0}}-e^{{z_0}+z_2}+2 e^{{z_1}+{z_0}}-e^{z_2+{z_1}+{z_0}}+e^{{z_1}+2{z_0}}+
 e^{z_2+z_1+2 {z_0}}\right).
\end{eqnarray}
Now, using (\ref{eq:poisson1}), and observing that
$$
P Q =\left(\begin{array}{cc}e^{-\frac{{z_2}}{2}+\frac{{z_1}}{2}}+e^{-\frac{{z_2}}{2}-\frac{{z_1}}{2}-{z_0}}+
e^{-\frac{{z_2}}{2}+\frac{{z_1}}{2}-{z_0}}
&-e^{-\frac{{z_2}}{2}+\frac{{z_1}}{2}}-e^{-\frac{{z_2}}{2}+\frac{{z_1}}{2}+{z_0}}\\
 -e^{-\frac{{z_2}}{2}+\frac{{z_1}}{2}}-e^{\frac{{z_2}}{2}+
 \frac{{z_1}}{2}}-e^{-\frac{{z_2}}{2}-\frac{{z_1}}{2}-z}-e^{-\frac{{z_2}}{2}+\frac{{z_1}}{2}-{z_0}}&
 e^{\frac{-{z_2}}{2}+\frac{{z_1}}{2}}+e^{\frac{-{z_2}}{2}+\frac{{z_1}}{2}+{z_0}}+
 e^{\frac{{z_2}}{2}+\frac{{z_1}}{2}+{z_0}}\\
\end{array}\right),
$$
$$
P Q^{-1}=\left(\begin{array}{cc}0&-e^{-\frac{{z_2}}{2}-\frac{{z_1}}{2}}\\
 e^{\frac{{z_2}}{2}+\frac{{z_1}}{2}}&e^{-\frac{{z_2}}{2}-\frac{{z_1}}{2}}+e^{\frac{{z_2}}{2}-\frac{{z_1}}{2}}+e^{\frac{{z_2}}{2}+\frac{{z_1}}{2}}\\
\end{array}\right),
$$
we obtain the same result.
\end{example}

Let us show that this Poisson bracket is skew symmetric: in fact
graphically $\{G_Q,G_P\} $ means that $\gamma_Q$ is now on top: fix
any direction on $\gamma_Q$ and  assign a $+$ when we turn left and
a minus when we turn right as above.
\begin{figure}[h]\begin{center}
\begin{pspicture}(-3.3,-4)(0.98,1)
\rput(-3.2,-1.5){$\Bigg\{$}
\psellipse[linecolor=red,linewidth=1pt](-2,-2)(-1,2)
 \pscircle[fillstyle=solid,fillcolor=white,linecolor=white](-1.1,-2.9){0.2}
  \psellipse[linecolor=blue,linewidth=1pt](-1,-1)(1,2)
   \rput(-0.2,-3.0){$\textcolor{blue}{G_Q}$}
\rput(-2.8,-3.8){$\textcolor{red}{G_P}$}
 \pscircle[linecolor=black,linewidth=2pt](-2,-0.2){0.7}\rput(0.5,-1.5){$\Bigg\}=\,\frac{1}{2}$}
\end{pspicture}
\begin{pspicture}(-3,-4)(0.8,1)
  \psellipse[linecolor=blue,linewidth=1pt](-1,-1)(1,2)
 \pscircle[fillstyle=solid,fillcolor=white,linecolor=white](-1.1,-2.9){0.2}
\psellipse[linecolor=red,linewidth=1pt](-2,-2)(-1,2)
 \pscircle[fillstyle=solid,fillcolor=white,linecolor=black,linewidth=2pt](-2,-0.2){0.7}
\psarc[linewidth=1.5pt](-1.3,-1.22){0.75}{96}{154}
\psarc[linewidth=1.5pt](-2.6,.4){0.9}{268}{0}
\rput(-0.6,-3.5){$G_{PQ^{-1}}$}
\rput(0.5,-1.5){$-\frac{1}{2}$}
\end{pspicture}
\begin{pspicture}(-3,-4)(0,1)
  \psellipse[linecolor=blue,linewidth=1pt](-1,-1)(1,2)
 \pscircle[fillstyle=solid,fillcolor=white,linecolor=white](-1.1,-2.9){0.2}
\psellipse[linecolor=red,linewidth=1pt](-2,-2)(-1,2)
 \pscircle[fillstyle=solid,fillcolor=white,linecolor=black,linewidth=2pt](-2,-0.2){0.7}
 \psarc[linewidth=1.5pt](-1,.2){0.75}{166}{240}
 \psarc[linewidth=1.5pt](-2.6,-1.2){0.75}{30}{90}
 \rput(-0.6,-3.5){$G_{PQ}$}
\end{pspicture}
\caption{Poisson bracket.}\label{pic:poisson1}
\end{center}
\end{figure}

The main problem is to close this Poisson algebra by using the skein
relation. That is, we want to choose a ``basis" of geodesic
functions, such that their Poisson brackets are expressed in terms
of the basis elements only. We quote the word ``basis"  because the
number of elements in it is actually larger than the dimension.

\begin{example} In our torus with one hole we choose the "basis"  given by $G_{z_2 z_1}$,
$G_{z_1 z_0}$,  $G_{z_0z_2}$. We get
$$
\{G_{z_2 z_1},G_{z_1 z_0}\} = \frac{1}{2} G_{z_0z_2} - \frac{1}{2} G_{z_2 z_1 z_0 z_1},
$$
and by using the skein relation
$$
\{G_{z_2 z_1},G_{z_1 z_0}\} = \frac{1}{2} G_{z_2 z_1}G_{z_1 z_0} - \frac{1}{2} G_{z_0z_2},
$$
The central element of this algebra is
$$
G_{z_2 z_1}^2 +G_{z_1 z_0}^2 + G_{z_0z_2}^2-G_{z_2 z_1}G_{z_1 z_0} G_{z_0z_2}.
$$
\end{example}

\begin{definition}
A geodesic is called {\it graph simple}\/ if it does not pass twice through the same edge.
\end{definition}

\begin{example}
In our torus with one hole we have only $3$ graph-simple geodesics,
the $A$-cycle $A$, the $B$-cycle $B$  and $AB^{-1}$. Their lengths
generate the algebra
\begin{eqnarray}\nn
&&
\{G_A,G_B\} = \frac{1}{2} G_A G_B- G_{A B^{-1}},\nn\\
&&
\{G_B,G_{AB^{-1}}\} = \frac{1}{2} G_B G_{A B^{-1}}-G_A,\nn\\
&&
\{G_{AB^{-1}},G_A\} = \frac{1}{2} G_A G_{A B^{-1}}-G_B.\nn
\end{eqnarray}
\end{example}

In the case of higher genus and $s=1,2$, Chekhov and Fock have
proved that the algebra is always closed by choosing the set of
graph--simple geodesics.

\begin{definition}
A {\it geodesic lamination}\/ $GL$ is a set of
non--(self)--intersecting geodesics.  The associated algebraic
object $\mathcal{GL}$ is the product of their lengths.
$$
\mathcal{GL}:=\prod_{\gamma\in GL} G_\gamma.
$$
\end{definition}

The advantage of considering geodesic laminations is that all relations are linear:
$$
\GL_1\GL_2=\sum_j c_{12}^j \GL_j,\qquad
\{\GL_1,\GL_2\}=\sum_j \tilde c_{12}^j \GL_j,
$$
where the constants $c^j_{il}\in{\mathbb Z}$ and  $\tilde c^j_{il}\in\frac{1}{2}{\mathbb Z}$.

The disadvantage is that the geodesic laminations can be chosen is infinite ways, therefore hey are not algebraically independent.

\section{Quantization}
 We quantize using the correspondence principle: to each coordinate
$z_\alpha$ we associate a Hermitean operator $Z^\hbar_\alpha$ such
that  the set of these Hermitean operators forms a
$C^\star$--algebra (in our case the algebra of operators on
$L^2(3g-3+s)$). We postulate the commutation relations between
Hermitean operators as follows:
\begin{equation}\label{comm-rel}
[Z_\alpha^\hbar,Z_\beta^\hbar]=2\pi\,i\hbar\{z_\alpha,z_\beta\}.
\end{equation}
The central elements of this algebra are the quantum analogues of the
Casimirs. They are mapped to
$$
P_j^\hbar=\sum_{i\in I_j} Z^\hbar_i,\qquad j=1,\dots s,
$$
where $I_j$ is the subset of indices labeling all edges in the face
containing the $j$-th hole. Observe that in the quantum picture we
loose the geometric picture: we're left only with an abstract
algebra.

The quantum morphisms work in the same way as before, with the
quantum analogue of $\phi$ which we denote by $\Phi^\hbar$. This
function must satisfy the following properties:
\begin{enumerate}
\item $\Phi^\hbar(z)-\Phi^\hbar(-z)=z$ to preserve the commutation relations.
\item Semiclassical limit: $\lim_{\hbar\to0}\Phi^\hbar(z)=\varphi(z)$.
\item Quantum pentagon identity. This is non-trivial as the quantized variables $Z^\hbar_\alpha$ don't commute.
\end{enumerate}

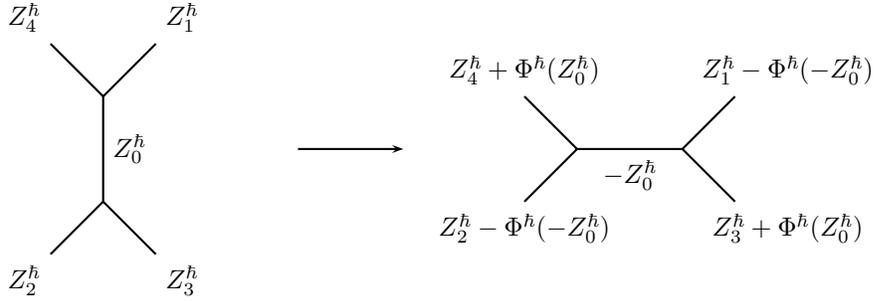
\begin{figure}[h]
\begin{center}
\psset{unit=7mm}
\begin{pspicture}(-1.5,-3.5)(13,1.3)
 \psline(-1,-3)(0,-2)(1,-3) \rput(-1.5,-3.5){$Z^\hbar_2$} \rput(1.5,-3.5){$Z^\hbar_3$}
 \psline(0,-2)(0,0) \rput(0.5,-1){$Z^\hbar_0$}
 \psline(-1,1)(0,0)(1,1) \rput(-1.5,1.5){$Z^\hbar_4$} \rput(1.5,1.5){$Z^\hbar_1$}
 \psline{->}(3.7,-1)(5.7,-1)
 \psline(8,-2)(9,-1)(8,0) \rput(8,0.5){$Z^\hbar_4+\Phi^\hbar(Z^\hbar_0)$}
 \rput(8,-2.5){$Z^\hbar_2-\Phi^\hbar(-Z^\hbar_0)$}
 \psline(9,-1)(11,-1) \rput(10,-1.5){$-Z^\hbar_0 $}
 \psline(12,-2)(11,-1)(12,0) \rput(13,0.5){$Z^\hbar_1-\Phi^\hbar(-Z^\hbar_0)$}
 \rput(13,-2.5){$Z^\hbar_3+\Phi^\hbar(Z^\hbar_0)$}
\end{pspicture}\caption{The quantum  flip in the general situation.}\label{q-flip}
\end{center}
\end{figure}

Let us see how to prove the first property.  We assume the quantum
flip to work similarly to the classical one, as described in figure
\ref{q-flip}.

Say we wish to preserve the commutation relation
$$
[Z_1^\hbar,Z_4^\hbar]=2\pi\, i.
$$
As we have seen just before Theorem \ref{th:morph},  $Z_1^\hbar\to Z_4^\hbar+\Phi^\hbar(Z_0^\hbar)$
and $Z_4^\hbar\to Z_2^\hbar-\Phi^\hbar(-Z_0^\hbar)$, so that we want
$$
[Z_4^\hbar+\Phi^\hbar(Z_0^\hbar),Z_2^\hbar-\Phi^\hbar(-Z_0^\hbar)]=2\pi\, i.
$$
Since $[Z^\hbar_4,Z_2^\hbar]=0$, this means that we want:
$$
[Z_4^\hbar,-\Phi^\hbar(-Z_0^\hbar)]+[\Phi^\hbar(Z_0^\hbar),Z_2^\hbar]=2\pi\, i,
$$
and by the postulated commutation relations (\ref{comm-rel}), we get
$$
\frac{d}{d Z_0^\hbar}\left(-\left(\Phi^\hbar(-Z^\hbar_0)\right)_{} +\left(\Phi^\hbar(Z^\hbar_0)\right)_{}\right)  =1,
$$
which proves the first property after integration in $Z_0^\hbar$.

One can prove that a good candidate for $\Phi^\hbar$ is the {\it quantum dilogarithm function} \cite{fad} (a.k.a.
the Barnes function in mathematical literature):
$$
\Phi^\hbar(z):= -\frac{\pi \hbar}{2}\int \frac{e^{-i p z}{\rm d}p}{\sinh(\pi p)\sinh(\pi p \hbar)}.
$$
This function has two properties of quasi--periodicity:
$$
\Phi^\hbar(z+i\pi\hbar)- \Phi^\hbar(z-i\pi\hbar)=\frac{2\pi i\hbar}{1+e^{-z}},\qquad
\Phi^\hbar(z+i\pi)- \Phi^\hbar(z-i\pi)=\frac{2\pi i}{1+e^{-\frac{z}{\hbar}}}.
$$

The quantum geodesic functions will be operators. Recall that at
classical level the geodesic functions were
$G\gamma=2\cosh(\frac{l_\gamma}{2})={\rm Tr}(P_{z_1,\dots,z_n})$.
When we quantize the quantized $Z^\hbar_\alpha$ do not commute
anymore: we have to choose an ordering in the trace. We'll denote
the quantum ordering by
$$
\begin{array}{c}{\rm x}\\ {\rm x}\\
\end{array} {\rm Tr} (P_{Z_1\dots Z_n})\begin{array}{c}{\rm x}\\ {\rm x}\\
\end{array}
$$
In the classical case we had that $G\gamma$ was a sum of exponential
terms of integer or half integer linear combinations of
$z_1,\dots,z_n$. Now in the quantum case we can no longer use the
property that $e^a e^b=e^{a+b}$, so we'll need to add a term $2\pi i
\hbar C(\gamma,\alpha)$, which will depend on the geodesic $\gamma$
and on the set of edges $\alpha$ it passes trough.

This is difficult to compute. However thanks to the fact that the
quantized operators $Z^\hbar_\alpha$ are Hermitean, one can prove a
few results. In particular, if we restrict to the class of  {\it
graph--simple geodesics ,}\/ i.e. geodesics which pass trough one
edge only once, then $C(\gamma,\alpha)=0$ so that the quantum
ordering is Weyl ordering and we have the Moyal product.

\begin{figure}[h]\begin{center}
\begin{pspicture}(-3.3,-4)(1.4,1)
\rput(-3.2,-1.5){$\Bigg[$}
  \psellipse[linecolor=blue,linewidth=1pt](-1,-1)(1,2)
 \pscircle[fillstyle=solid,fillcolor=white,linecolor=white](-1.1,-2.9){0.2}
  \rput(-0.2,-3.0){$\textcolor{blue}{G_Q}$}
\psellipse[linecolor=red,linewidth=1pt](-2,-2)(-1,2)
\rput(-2.8,-3.8){$\textcolor{red}{G_P}$}
 \pscircle[linecolor=black,linewidth=2pt](-2,-0.2){0.7}\rput(0.7,-1.5){$\Bigg]=\,q^{-\frac{1}{2}}$}
\end{pspicture}
\begin{pspicture}(-3,-4)(0.8,1)
  \psellipse[linecolor=blue,linewidth=1pt](-1,-1)(1,2)
 \pscircle[fillstyle=solid,fillcolor=white,linecolor=white](-1.1,-2.9){0.2}
\psellipse[linecolor=red,linewidth=1pt](-2,-2)(-1,2)
 \pscircle[fillstyle=solid,fillcolor=white,linecolor=black,linewidth=2pt](-2,-0.2){0.7}
 \psarc[linewidth=1.5pt](-1,.2){0.75}{166}{240}
 \psarc[linewidth=1.5pt](-2.6,-1.2){0.75}{30}{90}
 \rput(-0.6,-3.5){$G_{PQ}$}
\rput(0.5,-1.5){$+\,q^{\frac{1}{2}}$}
\end{pspicture}
\begin{pspicture}(-3,-4)(0,1)
  \psellipse[linecolor=blue,linewidth=1pt](-1,-1)(1,2)
 \pscircle[fillstyle=solid,fillcolor=white,linecolor=white](-1.1,-2.9){0.2}
\psellipse[linecolor=red,linewidth=1pt](-2,-2)(-1,2)
 \pscircle[fillstyle=solid,fillcolor=white,linecolor=black,linewidth=2pt](-2,-0.2){0.7}
\psarc[linewidth=1.5pt](-1.3,-1.22){0.75}{96}{154}
\psarc[linewidth=1.5pt](-2.6,.4){0.9}{268}{0}
\rput(-0.6,-3.5){$G_{PQ^{-1}}$}
\end{pspicture}
\caption{Quantum version of skein relation and Poisson bracket.}\label{pic:q-poisson}
\end{center}
\end{figure}
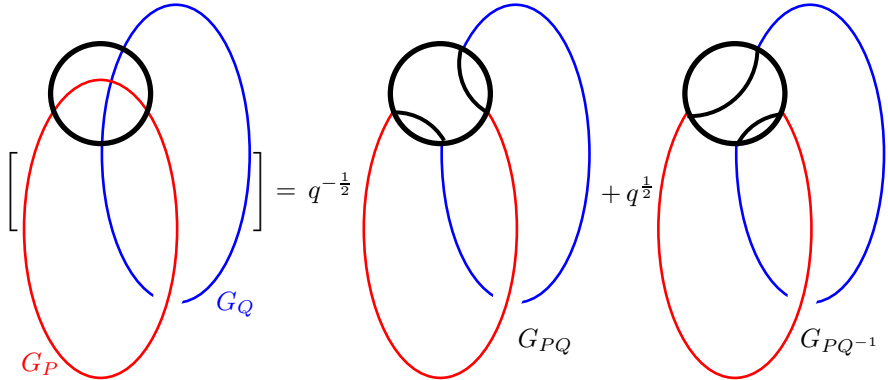

Skein relation and Poisson bracket are now encoded in a single
{\it quantum skein relation,} which for one crossing is given in
figure \ref{pic:q-poisson} in which  $q=e^{-i\pi\hbar}$, so that for
$q=1$ we get the classical skein relation and the first correction in
$\hbar$ produces the Poisson bracket.

\begin{figure}[h]\begin{center}
\begin{pspicture}(-1.7,-2)(2.6,2)
\rput(-2.7,1){$\Bigg[$}
\psline[linecolor=blue,linewidth=1pt](-2.5,1)(2.5,1)
\pscircle[fillstyle=solid,fillcolor=white,linecolor=white](-1.75,1){0.1}
\pscircle[fillstyle=solid,fillcolor=white,linecolor=white](1.75,1){0.1}
\psarc[linecolor=red,linewidth=1pt](-0,0){2}{0}{180}
\rput(2.7,1){$\Bigg]$}
\rput(3.2,1){$=$}
\end{pspicture}
\begin{pspicture}(-3.4,-2)(1.6,2)
\psline[linecolor=blue,linewidth=1pt](-2.5,1)(2.5,1)
\psarc[linecolor=red,linewidth=1pt](-0,0){2}{0}{180}
\pscircle[fillstyle=solid,fillcolor=white,linecolor=white](-1.75,1){0.3}
\pscircle[fillstyle=solid,fillcolor=white,linecolor=white](1.75,1){0.3}
\psarc[linecolor=black,linewidth=1pt](-1.5,1.1){0.15}{112}{315}
\psarc[linecolor=black,linewidth=1pt](-2,.83){0.17}{-42}{110}
\psarc[linecolor=black,linewidth=1pt](2,.83){0.17}{76}{220}
\psarc[linecolor=black,linewidth=1pt](1.5,1.13){0.13}{-110}{60}
\rput(3,1){$+$}
\end{pspicture}
\begin{pspicture}(-1.7,-4)(2.6,-2)
\rput(-3,-2){$+\, q$}
\psline[linecolor=blue,linewidth=1pt](-2.5,-2)(2.5,-2)
\psarc[linecolor=red,linewidth=1pt](-0,-3){2}{0}{180}
\pscircle[fillstyle=solid,fillcolor=white,linecolor=white](-1.75,-2){0.3}
\pscircle[fillstyle=solid,fillcolor=white,linecolor=white](1.75,-2){0.3}
\psarc[linecolor=black,linewidth=1pt](-1.5,-1.9){0.15}{112}{315}
\psarc[linecolor=black,linewidth=1pt](-2,-2.17){0.17}{-42}{110}
\psarc[linecolor=black,linewidth=1pt](2,-1.5){0.5}{210}{276}
\psarc[linecolor=black,linewidth=1pt](1.36,-2.6){0.6}{26}{96}
\rput(3,-2){$+\, q^{-1}$}
\end{pspicture}
\begin{pspicture}(-3.4,-4)(1.6,-2)
\psline[linecolor=blue,linewidth=1pt](-2.5,-2)(2.5,-2)
\psarc[linecolor=red,linewidth=1pt](-0,-3){2}{0}{180}
\pscircle[fillstyle=solid,fillcolor=white,linecolor=white](-1.75,-2){0.3}
\pscircle[fillstyle=solid,fillcolor=white,linecolor=white](1.75,-2){0.3}
\psarc[linecolor=black,linewidth=1pt](-1.4,-2.49){0.5}{97}{160}
\psarc[linecolor=black,linewidth=1pt](-2,-1.5){0.5}{-96}{-30}
\psarc[linecolor=black,linewidth=1pt](2,-2.17){0.17}{76}{220}
\psarc[linecolor=black,linewidth=1pt](1.5,-1.87){0.13}{-110}{60}
\rput(3,-2){$+$}
\end{pspicture}
\begin{pspicture}(-1.5,-6)(2.6,-4)
\rput(-3,-5){$+$}
\psline[linecolor=blue,linewidth=1pt](-2.5,-5)(2.5,-5)
\psarc[linecolor=red,linewidth=1pt](-0,-6){2}{0}{180}
\pscircle[fillstyle=solid,fillcolor=white,linecolor=white](-1.75,-5){0.3}
\pscircle[fillstyle=solid,fillcolor=white,linecolor=white](1.75,-5){0.3}
\psarc[linecolor=black,linewidth=1pt](-1.4,-5.49){0.5}{97}{160}
\psarc[linecolor=black,linewidth=1pt](-2,-4.5){0.5}{-96}{-30}
\psarc[linecolor=black,linewidth=1pt](2,-4.5){0.5}{210}{276}
\psarc[linecolor=black,linewidth=1pt](1.36,-5.6){0.6}{26}{96}
\rput(3,-5){$=$}
\end{pspicture}
\begin{pspicture}(-3.4,-2)(1.6,-0)
\psline[linecolor=blue,linewidth=1pt](-2.5,0)(2.5,0)
\psarc[linecolor=red,linewidth=1pt](-0,-2.3){2}{0}{180}
\end{pspicture}

\caption{Quantum Whitehead move I: it holds if the empty loop is $-q-q^{-1}$.}\label{empty-loop}
\end{center}
\end{figure}
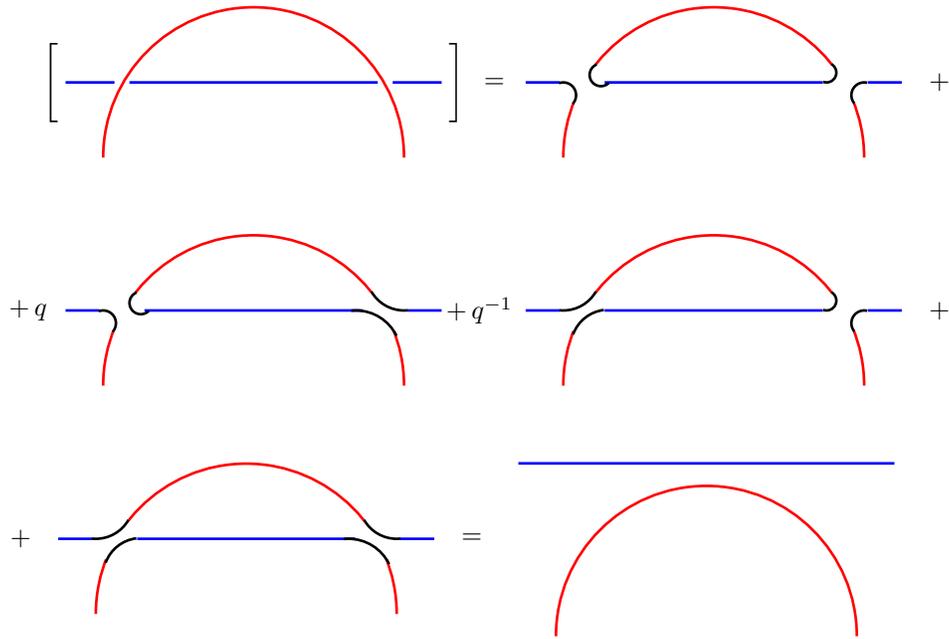

We need to assign a geodesic length to the empty loop too. This is
fixed to be equal to $-2$ in the classical case, and it becomes
$-q-q^{-1}$ in the quantum case. The reason for this is described in figure
\ref{empty-loop}, in which we consider two a~priori commuting
geodesic functions and compute their quantum skein relation.

Let us see how to close the algebra in the quantum case. In fact
before we had two relations, i.e., the Poisson bracket and the skein
relation. When computing the former, we used the latter to replace
unwanted terms. But now we have only one relation, so what shall we
do? The trick is to use the non-commutativity, i.e. to use the
quantum skein relation twice, once for $G_A G_B$ and once for $G_B
G_A$:
$$
G_A G_B = q^{-\frac{1}{2}} G_{AB} + q^{\frac{1}{2}}  G_{A B^{-1}},\qquad
G_B G_A = q^{\frac{1}{2}} G_{AB} + q^{-\frac{1}{2}}  G_{A B^{-1}},
$$
so that if we define the quantum commutator as
$$
[G_A,G_B]_q:=  q^{\frac{1}{2}} G_{A} G_B - q^{-\frac{1}{2}}  G_B G_A,
$$
we get
$$
[G_A,G_B]_q:=  (q^{\frac{1}{2}} - q^{-\frac{1}{2}})G_{A B^{-1}}.
$$


\def\thetheorem{A.\arabic{theorem}}
\def\theprop{A.\arabic{prop}}
\def\thelemma{A.\arabic{lemma}}
\def\thecor{A.\arabic{cor}}
\def\theexam{A.\arabic{exam}}
\def\theremark{A.\arabic{remark}}
\def\theequation{A.\arabic{equation}}

\section{Bibliographic references}

The combinatorial description of Teichm\"uller theory was proposed by R.C. Penner
\cite{Penn1,pen1,pen2} in the case of punctured surfaces and was advanced by V.V. Fock
\cite{Fock93,Fock97} for the case of surfaces with holes.

The Poisson algebras of geodesic functions have a long history, first, irrespectively on
the combinatorial description. Their investigation has been started in the $2+1$ gravity pattern
in which W. Goldman \cite{Gold1} established his famous Goldman bracket. Using it, J. Nelson and
T. Regge (subsequently in collaboration with S. Carlip and Zertuche) found \cite{NR1,NR2} algebras
for a special basis of geodesic functions. The geodesic algebras were found in the graph
(combinatorial) description in \cite{CF2}, whereas the representation in terms of the graph-simple
geodesics for a general Nelson--Regge algebra was found in \cite{CF3}.

The quantization of the Teichm\"uller space coordinates satisfying the pentagon identity was proposed
by Chekhov and Fock \cite{CF1} and independently by Kashaev \cite{Kash1}. The quantum dilogarithm function
in the modern concept was introduced by L.D. Faddeev \cite{fad} in his studies of the quantum Heisenberg
algebra symmetries. A good accounting for representations of quantum geodesic functions and of their
spectral properties is in \cite{Kash2}.

Further development of this activity was the quantization of the Thurston theory \cite{CP}, the Fock--Goncharov
construction for quantizing $PSL(n, \mathbb R)$ and other algebras \cite{FG}, generalizations to the
bordered surfaces case, etc. The Nelson--Regge algebras has appeared as algebras of the upper-triangular groupoid
\cite{Bondal} and of the Stokes parameters \cite{Ug} for the Frobenius manifolds.

\appendix

\section{The hyperbolic metric}\label{app-a}

The {\it hyperbolic distance}\/ $\rho(z,w)$ between two points
$z,w\in\mathbb H$  is defined as the length of the geodesic
connecting them. This defines $\mathbb H$ as a metric space. It is a
straightforward computation to show that
$$
\rho(z,w)= \ln\left(\frac{|z-\bar w|+|z-w|}{|z-\bar w|-|z-w|}
\right),
$$
so that
$$
\sinh^2\left(\frac{1}{2} \rho(z,w)
\right)= \frac{|z-w|^2}{4\, \IM(z)\,\IM(w)}.
$$
Using this formula we are now going to show that the topology
induced by the hyperbolic distance on $\mathbb H$ is the same as the
one induced by the Euclidean metric. In fact a hyperbolic circle is
also a Euclidean circle and vice versa. Let us see this in an
example: consider the hyperbolic circle $C$ of center $i$ and radius
$\delta$:
$$
C=\left\{z\in\mathbb H\big|\, \sinh^2\left(\frac{1}{2} \rho(z,i)
\right)=\sinh^2\left(\frac{\delta}{2} \right)\right \}.
$$

\begin{figure}[h]
\begin{center}
\begin{pspicture}(-2,-1)(12,4)
 \psline{->}(2,-1)(6,-1)
 \psline{->}(3.5,-1)(3.5,3.5)
\pscircle(3.5,1){1.5}
\psdot(3.5, 0)\rput(3.7,0){$i$}
\psdot(3.5,1)\rput(4.2,1){$i \cosh\delta$}
\end{pspicture}
\end{center}
 \caption{}\label{pic:poincare}
\end{figure}

Clearly we get
\begin{eqnarray}\nn
&&
C= \left\{z=x+i y\, \big|\, |z-i|^2 = 4 y \sinh^2\left(\frac{\delta}{2} \right)\right \}=\nn\\
&&
\quad = \left\{z=x+i y\, \big|\, x^2+y^2+1 =
2 y \left(2  \sinh^2\left(\frac{\delta}{2} \right) +1\right)= 2 y \cosh \delta \right\}=\nn\\
&&
\quad =\left\{z=x+i y\, \big|\,  x^2+(y-\cosh \delta)^2=\cosh^2\delta-1=\sinh^2\delta\right\}.
\nn
\end{eqnarray}
Therefore $C$ is a Euclidean circle of center $z=i \cosh\delta$ and radius $\sinh\delta$.

\section{Poisson brackets between geodesic functions}

We now calculate the Poisson bracket between geodesic functions labeled ``1'' and ``2'' and having
a crossing as indicated in figure \ref{pic:crossing2}.

We use the following notation standard, for example, in integrable
systems. We treat matrix elements of the first and second
geodesic functions as linear operators in two different
two-dimensional spaces. The superscripts
``(1)'' and ``(2)'' indicate quantities belonging to the space with
the corresponding number in the direct product of spaces,
$e^{(s)}_{ij}$ is the matrix in the $s$th space $(s=1,2)$ with the
only nonzero entry to be the unity on the crossing of $i$th row and
$j$th column, and $\varepsilon_{ij}$ is the totally antisymmetric
tensor (with $\varepsilon_{12}=1$. Every operator acting in the
direct product of spaces can be decomposed over the basis with the
basic vectors to be $e^{(1)}_{ij}\otimes e^{(2)}_{kl}$, $i,j,k,l=1,2$.
We also distinguish between
evaluating traces in these two spaces.

\begin{figure}[h]
\begin{center}
\begin{pspicture}(-7,-2)(6,2.5)
\psline[linewidth=1pt](-3,2)(-1.5,0.5)
\psline[linewidth=1pt](-3,-2)(-1.5,-0.5)
\psline[linewidth=1pt](3,2)(1.5,0.5)
\psline[linewidth=1pt](3,-2)(1.5,-0.5)
\psline[linewidth=1pt](-3.5,1.5)(-2,0)
\psline[linewidth=1pt](-3.5,-1.5)(-2,0)
\psline[linewidth=1pt](3.5,1.5)(2,0)
\psline[linewidth=1pt](3.5,-1.5)(2,0)
\psline[linewidth=1pt](-1.5,.5)(1.5,.5)
\psline[linewidth=1pt](-1.5,-.5)(1.5,-.5)
\psbezier[linecolor=blue, linewidth=1pt](-3.3,-1.7)(-1.8,-0.2)(-1.2,0)(0,0)
\psbezier[linecolor=blue, linewidth=1pt](0,0)(1.2,0)(1.8,0.2)(3.3,1.7)
\psbezier[linecolor=white, linewidth=4pt](-3.3,1.7)(-1.8,0.2)(-1.2,0)(0,0)
\psbezier[linecolor=white, linewidth=4pt](0,0)(1.2,0)(1.8,-0.2)(3.3,-1.7)
\psbezier[linecolor=red, linewidth=1pt](-3.3,1.7)(-1.8,0.2)(-1.2,0)(0,0)
\psbezier[linecolor=red, linewidth=1pt](0,0)(1.2,0)(1.8,-0.2)(3.3,-1.7)
\rput(-2.5,2){\makebox(0,0){$z_1$}}
\rput(2.5,2){\makebox(0,0){$z_2$}}
\rput(2.5,-2){\makebox(0,0){$z_3$}}
\rput(-2.5,-2){\makebox(0,0){$z_4$}}
\rput(0,0.8){\makebox(0,0){$z_5$}}
\rput(-3.5,2){\makebox(0,0){${\mathbf 1}$}}
\rput(-3.5,-2){\makebox(0,0){${\mathbf 2}$}}
\end{pspicture}
\end{center}
 \caption{}\label{pic:crossing2}
\end{figure}

Another way to represent $X_z$ is
$$
X_z^{(1)}=\sum_{i,j=1}^2(-1)^{\varepsilon_{ij}+1}e^{\varepsilon_{ij}z/2}e^{(1)}_{ij}.
$$
Then,
$$
\{X_{z_1}^{(1)},X_{z_4}^{(2)}\}=\frac{-1}{4}X_{z_1}^{(1)}\otimes X_{z_4}^{(2)}{\mathcal S}^{12},
$$
where $z_1$ and $z_4$ are as in figure \ref{pic:crossing2} and
$$
{\mathcal S}^{12}:=e^{(1)}_{11}\otimes e^{(2)}_{11}-e^{(1)}_{22}\otimes e^{(2)}_{11}
-e^{(1)}_{11}\otimes e^{(2)}_{22}+e^{(1)}_{22}\otimes e^{(2)}_{22}
$$
is a diagonal matrix. Note that
\begin{equation}
X_{z_1}^{(1)}\otimes X_{z_4}^{(2)}{\mathcal S}^{12}={\mathcal S}^{12}X_{z_1}^{(1)}\otimes X_{z_4}^{(2)},
\label{XXS}
\end{equation}
so we can put the matrix ${\mathcal S}^{12}$ from any side of this tensorial product.
When calculating the Poisson brackets between geodesic functions $G_1=\tr P$ and $G_2=\tr Q$ in figure \ref{pic:crossing2}, we
insert a proper number (six in this case) of the matrices ${\mathcal S}^{12}$ in the tensorial product of
$$
P^{(1)}\otimes Q^{(2)}=\cdots X_{z_3}^{(1)}R^{(1)}X_{z_5}^{(1)}L^{(1)}X_{z_1}^{(1)}\cdots \otimes
\cdots X_{z_2}^{(2)}L^{(2)}X_{z_5}^{(2)}R^{(2)}X_{z_4}^{(2)}\cdots;
$$
due to the property (\ref{XXS}) all these insertions can be done to the right from $X_{z_3}^{(1)}$ and $X_{z_2}^{(2)}$ and
to the left from $X_{z_1}^{(1)}$ and $X_{z_4}^{(2)}$. We show these insertions schematically on the figure \ref{pic:XXS},
where the signs $(+)$ and $(-)$ indicate with which sign the corresponding matrix ${\mathcal S}^{12}$ enters the direct
product. Eventually, using that $X_zX_z=RL=LR=-{\mathbf 1}$, we push all ${\mathcal S}^{12}$ insertions to the single place
just to the right from $X_{z_5}^{(1)}\otimes X_{z_5}^{(2)}$. We then have for the bracket
\begin{eqnarray}
\nonumber
&{}&\cdots \frac14 (RX_{z_5})^{(1)}\otimes (LX_{z_5})^{(2)}
\biggl\{(X_{z_5}L)^{(1)}\otimes {\mathbf 1}^{(2)}{\mathcal S}^{12}(RX_{z_5})^{(1)}\otimes {\mathbf 1}^{(2)}\biggr.
\\
\nonumber
&{}&-(X_{z_5}L)^{(1)}\otimes (X_{z_5}R)^{(2)}{\mathcal S}^{12}(RX_{z_5})^{(1)}\otimes (LX_{z_5})^{(2)}
+X_{z_5}^{(1)}\otimes R^{(2)}{\mathcal S}^{12}X_{z_5}^{(1)}\otimes L^{(2)}
\\
\nonumber
&{}&+{\mathbf 1}^{(1)}\otimes (X_{z_5}R)^{(2)}{\mathcal S}^{12}{\mathbf 1}^{(1)}\otimes (LX_{z_5})^{(2)}
+L^{(1)}\otimes X_{z_5}^{(2)}{\mathcal S}^{12}R^{(1)}\otimes X_{z_5}^{(2)}
\\
&{}&\biggl.
-L^{(1)}\otimes R^{(2)}{\mathcal S}^{12}R^{(1)}\otimes L^{(2)}\biggr\}L^{(1)}\otimes R^{(2)}\cdots\, .
\nonumber
\end{eqnarray}

\begin{figure}[h]
\begin{center}
\begin{pspicture}(-7,-2)(6,2.5)
\psline[linestyle=dashed, linewidth=1pt](-2.5,1.7)(-2.5,-1.7)
\psline[linestyle=dashed, linewidth=1pt](2.5,1.7)(2.5,-1.7)
\psbezier[linestyle=dashed, linewidth=1pt](-2.2,1.7)(-2,0.2)(0.5,-0.2)(0.7,-1.7)
\psbezier[linestyle=dashed, linewidth=1pt](-2.2,-1.7)(-2,-0.2)(0.5,0.2)(0.7,1.7)
\psbezier[linestyle=dashed, linewidth=1pt](2.2,1.7)(2,0.2)(-0.5,-0.2)(-0.7,-1.7)
\psbezier[linestyle=dashed, linewidth=1pt](2.2,-1.7)(2,-0.2)(-0.5,0.2)(-0.7,1.7)
\rput(-3,1.5){\makebox(0,0){$X_{z_3}$}}
\rput(-1.5,1.5){\makebox(0,0){$R$}}
\rput(0,1.5){\makebox(0,0){$X_{z_5}$}}
\rput(1.5,1.5){\makebox(0,0){$L$}}
\rput(3,1.5){\makebox(0,0){$X_{z_1}$}}
\rput(-3,-1.5){\makebox(0,0){$X_{z_2}$}}
\rput(-1.5,-1.5){\makebox(0,0){$L$}}
\rput(0,-1.5){\makebox(0,0){$X_{z_5}$}}
\rput(1.5,-1.5){\makebox(0,0){$R$}}
\rput(3,-1.5){\makebox(0,0){$X_{z_4}$}}
\rput(-3.8,1.5){\makebox(0,0){${\mathbf 1}$}}
\rput(-3.8,-1.5){\makebox(0,0){${\mathbf 2}$}}
\rput(-2.6,2){\makebox(0,0){$\scriptsize (+)$}}
\rput(-2.1,2){\makebox(0,0){$\scriptsize (-)$}}
\rput(-0.7,2){\makebox(0,0){$\scriptsize (+)$}}
\rput(0.7,2){\makebox(0,0){$\scriptsize (+)$}}
\rput(2.1,2){\makebox(0,0){$\scriptsize (+)$}}
\rput(2.6,2){\makebox(0,0){$\scriptsize (-)$}}
\end{pspicture}
\end{center}
 \caption{}\label{pic:XXS}
\end{figure}

It remains just to calculate explicitly the part in the brackets. It reads
\begin{eqnarray}
\nonumber
&{}&
\left(%
\begin{array}{cc}
  1 & -2e^{z_5} \\
  0 & -1 \\
\end{array}%
\right)
\otimes
\left(%
\begin{array}{cc}
  1 & 0 \\
  2e^{-z_5} & -1 \\
\end{array}%
\right)
+
\left(%
\begin{array}{cc}
  1 & -2e^{z_5} \\
  0 & -1 \\
\end{array}%
\right)
\otimes
\left(%
\begin{array}{cc}
  1 & 0 \\
  0 & -1 \\
\end{array}%
\right)
\\
&{}&
\nonumber
+
\left(%
\begin{array}{cc}
  1 & 0 \\
  0 & -1 \\
\end{array}%
\right)
\otimes
\left(%
\begin{array}{cc}
  1 & 2 \\
  0 & -1 \\
\end{array}%
\right)
+
\left(%
\begin{array}{cc}
  1 & 0 \\
  0 & -1 \\
\end{array}%
\right)
\otimes
\left(%
\begin{array}{cc}
  1 & 0 \\
  2e^{-z_5} & -1 \\
\end{array}%
\right)
\\
\nonumber
&{}&
+
\left(%
\begin{array}{cc}
  1 & 0 \\
  -2 & -1 \\
\end{array}%
\right)
\otimes
\left(%
\begin{array}{cc}
  1 & 0 \\
  0 & -1 \\
\end{array}%
\right)
-
\left(%
\begin{array}{cc}
  1 & 0 \\
  -2 & -1 \\
\end{array}%
\right)
\otimes
\left(%
\begin{array}{cc}
  1 & 2 \\
  0 & -1 \\
\end{array}%
\right)
\\
\nonumber
&=&
4\left(%
\begin{array}{cc}
  0 & e^{z_5} \\
  0 & 0 \\
\end{array}%
\right)
\otimes
\left(%
\begin{array}{cc}
  0 & 0 \\
  e^{-z_5} & 0 \\
\end{array}%
\right)
+
4\left(%
\begin{array}{cc}
  0 & 0 \\
  1 & 0 \\
\end{array}%
\right)
\otimes
\left(%
\begin{array}{cc}
  0 & 1 \\
  0 & 0 \\
\end{array}%
\right)
\\
\nonumber
&{}&\qquad
+
2\left(%
\begin{array}{cc}
  1 & 0 \\
  0 & -1 \\
\end{array}%
\right)
\otimes
\left(%
\begin{array}{cc}
  1 & 0 \\
  0 & -1 \\
\end{array}%
\right)
\\
\nonumber
&=&
4\bigl(e_{12}^{(1)}\otimes e_{21}^{(2)}+e_{21}^{(1)}\otimes e_{12}^{(2)}+e_{11}^{(1)}\otimes e_{11}^{(2)}+e_{22}^{(1)}\otimes e_{22}^{(2)}\bigr)
-2\bigl({\mathbf 1}^{(1)}\otimes {\mathbf 1}^{(2)}\bigr)
\\
&=&4{\mathcal R}^{12}-2\,{\mathbf 1}^{12},
\nonumber
\end{eqnarray}
where ${\mathcal R}^{12}$ (the expression in the first brackets) is
the famous permutation $R$-matrix\footnote{This matrix, permuting
the spaces $1$ and $2$, has the form
$\sum_{i,j=1}^2e^{(1)}_{ij}\otimes e^{(2)}_{ji}$, and it is the
simplest nontrivial matrix satisfying the celebrated Yang--Baxter
equation ${\mathcal R}^{12}{\mathcal R}^{23}{\mathcal
R}^{12}={\mathcal R}^{23}{\mathcal R}^{12}{\mathcal R}^{23}$ in the
direct product of three spaces.}

Taking into account the factor $1/4$, we immediately obtain that
\begin{eqnarray}
\nonumber
&{}&
\bigl\{\tr P^{(1)}, \tr Q^{(2)}\bigr\}=\tr_1\tr_2\bigl[ P^{(1)}\otimes Q^{(2)} \bigl( {\mathcal R}^{12}-\tfrac12 {\mathbf 1}^{12}\bigr)\bigr]
\\
\nonumber
&=&\tr PQ-\frac12 \tr P\cdot\tr Q,
\end{eqnarray}
with the usual traces in the last line. Using the skein relation
$$
\tr P\cdot \tr Q=\tr PQ+\tr PQ^{-1},
$$
we come to the standard form of the Poisson bracket:
$$
\{\tr P,\tr Q\}=\frac12 \tr PQ-\frac12 \tr PQ^{-1}.
$$


\def\refjl#1#2#3#4#5#6#7{\bibitem{#1}{\frenchspacing\rm#2}, #6,
{\frenchspacing\it#3\/}, {\bf#4} (#7) #5.}

\def\refbk#1#2#3#4#5{\bibitem{#1}{\frenchspacing\rm#2},
``{\sl#3\/}'', #4 (#5).}

\def\refpp#1#2#3#4#5{
\bibitem{#1}
\smallskip\noindent{\frenchspacing\rm#2}, {``#3''}\ #4\ (#5).}

\def\refcf#1#2#3#4#5#6#7{
\bibitem{#1}
\smallskip\noindent{\frenchspacing\rm#2}, {\rm#3},
in ``{\sl#4\/}'' [#5] #6\ (#7).}

\def\sideiii#1#2#3#4{\refcf{#1}{#2}{#3}{SIDE III --- Symmetries and
Integrability of Difference  Equations}{D.\ Levi and O.\ Ragnisco,
Editors}{{\frenchspacing\it CRM Proc. Lect. Notes Series\/}, vol.
{\bf25}, Amer. Math. Soc., Providence, RI, pp.\ #4}{2000}}

\vskip 2 mm \noindent{\bf Acknowledgments.}
We are grateful to James
Montaldi for his help with the pictures. This work was supported
by the Manchester Institute for Mathematical Sciences MIMMS, by the European Science Foundation Programme
``Methods of Integrable Systems, Geometry, Applied Mathematics"
(MISGAM), Marie Curie RTN ``European Network in Geometry,
Mathematical Physics and Applications"  (ENIGMA),  and by the EPSRC Fellowship EP/D071895/1.

\end{document}